\documentclass[a4paper, 11pt]{article}
\usepackage{amsmath,amsfonts,amssymb,amscd,amsthm,amsbsy}
\usepackage[dvips]{graphicx}
\graphicspath{{images/}}
\usepackage{mathabx} 
\usepackage{comment} 
\usepackage[dvips]{color}
\usepackage[dvipsnames]{xcolor} \usepackage{fancyhdr}
\usepackage{hyperref}
\usepackage{calrsfs}

\usepackage{listings}
\date{}

\newtheorem{theorem}{Theorem}[section]
\newtheorem{lemma}[theorem]{Lemma}

\theoremstyle{remark}
\newtheorem{remark}{Remark}[section]
\theoremstyle{definition}

\DeclareMathAlphabet{\pazocal}{OMS}{zplm}{m}{n}
\newcommand{\R}{\mathbb{R}}
\newcommand{\N}{\mathbb{N}}

\newcommand{\Z}{\mathbb{Z}}
\newcommand{\T}{\mathbb{T}}

\newcommand{\ep}{\varepsilon}
\newcommand{\kzero}{K_0}
\newcommand{\kh}{\widehat K}
\newcommand{\om}{\omega}
\newcommand{\omh}{\widehat \omega}
\newcommand{\kt}{\widetilde K}

\newcommand{\pt}{\Pi^{\top}}
\newcommand{\pp}{\Pi^{\bot}}
\newcommand{\gaep}{\Gamma^\ep}

\DeclareMathOperator{\Lip}{Lip}

\newcommand{\bydef}{\,\stackrel{\mbox{\tiny\textnormal{\raisebox{0ex}[0ex][0ex]{def}}}}{=}\,}

\begin{document}

\title
{Persistence of periodic orbits under state-dependent delayed
  perturbations: computer-assisted proofs}

\author{ %
  Joan Gimeno \footnote{University of Rome Tor Vergata, Dipartimento
  di Matematica, Via della Ricerca Scientifica, 00133, Rome,
  Italy. {\tt joan@maia.ub.es}}  
  \and %
  Jean-Philippe Lessard \footnote{McGill University, Department of
  Mathematics and Statistics, 805 Sherbrooke Street West, Montreal,
  QC, H3A 0B9, Canada. {\tt jp.lessard@mcgill.ca}}
  \and %
  J.D. Mireles James \footnote{Florida Atlantic University, Department
  of Mathematical Sciences, Science Building, Room 234, 777 Glades
  Road, Boca Raton, Florida, 33431 , USA. {\tt jmirelesjames@fau.edu}}
  \and %
  Jiaqi Yang \footnote{ICERM, Brown University, 121 South Main Street,
  Providence, RI, 02903, USA. {\tt jiaqi\_yang1@brown.edu}}
}

\maketitle

\begin{abstract}
A computer-assisted argument is given, which provides existence proofs for periodic orbits in state-dependent delayed 
perturbations of ordinary differential equations (ODEs). Assuming that the unperturbed ODE has an isolated periodic orbit, 
we introduce a set of polynomial inequalities whose successful verification leads to the existence of periodic orbits in the 
perturbed delay equation. We present a general algorithm, which describes a way of computing the coefficients of the 
polynomials and optimizing their variables so that the polynomial inequalities are satisfied. The algorithm uses the tools 
of validated numerics together with Chebyshev series expansion to obtain the periodic orbit of the ODE as well as the
 solution of the variational equations, which are both used to compute rigorously the coefficients of the polynomials. 
 We apply our algorithm to prove the existence of periodic orbits in a state-dependent delayed perturbation of the 
 van der Pol equation. 
\end{abstract}

\begin{center}
\end{center}

\begin{center}
\end{center}


\pagestyle{fancy}
\lhead{}
\chead{CAPs of Persistence of Periodic orbit}
\rhead{}


\section{Introduction} \label{sec:intro}
The goal of the present work is to develop computer assisted arguments for proving the existence
of periodic orbits in state-dependent delay differential equations (SDDEs)
arising as perturbations of ordinary differential equations (ODEs).
After fixing bounds on the norm of the perturbation, the delay, and their derivatives,  
our method determines (a) values of the perturbation parameter $\ep$ so that the periodic orbit of the ODE persists
into the SDDE, and (b) provides explicit bounds on the distance, 
in an appropriate norm, between the perturbed and  unperturbed solutions. 
Note that, while an ODE generates a finite 
dimensional dynamical system, an SDDE  --if it
generates a semi-flow at all-- has phase space as an
infinite dimensional Banach manifold 
\cite{MR2019242,MR4261206,MR3642771,MR1345150,MR2084753}.
This makes the perturbation arguments fairly delicate, and one novelty of 
our method is that it does not 
require to consider the Cauchy problem for the SDDE.

Regularity issues are also quite subtle.
For example, if the vector field generating the ODE is real analytic then 
any periodic orbits are real analytic as well.  Not so for SDDEs.  It has been shown 
that periodic solutions for SDDEs may be analytic in the neighborhood of
a certain point, and only $C^\infty$ on other portions of the orbit
 \cite{MR3229655,MR3748505,MR3992063}.  Moreover, it has been
 conjectured that this state of affairs is generic.  
Regularity questions have practical implications on the set-up of the 
perturbative argument -- for example we cannot 
employ analytic norms.

To manage these difficulties, we employ  
an approach based on the parameterization method
\cite{MR1976079,MR1976080,MR2177465,MR2240743,MR2289544,MR2299977}.
The parameterization method is a functional analytic framework for studying invariant manifolds,
which exploits the fact  that recurrent enough solutions often have much nicer properties
than solutions with arbitrary initial conditions.  The idea of the method is to 
formulate a chart or covering map 
for the invariant object as the solution of an invariance equation,
and studying such an equation allows the problem to be attacked using all the
tools of nonlinear analysis and computational mathematics.
A much more thorough description of the parameterization method 
with many applications is in
\cite{MR3467671}.

Recently a number of authors have made substantial progress using the 
parameterization method to study invariant manifolds in 
ill posed problems \cite{MR3900818,MR4066033,MR4120821,MR3749257},
and in particular the method has been used successfully to study
periodic and quasi-periodic solutions of SDDEs
and their attached stable/unstable manifolds
\cite{MR3501842,MR3736145,MR4112213,MR4287353,Per}.
We think of this  as an application of  the Poincar\'{e} program
in problems where the semi-flow theory is underdeveloped or 
otherwise problematic.  That is, one builds up an understanding of the 
dynamics one invariant object at a time, putting aside the fact that the 
dynamics in a full neighborhood of the
 invariant sets may or may not make sense at all.

The paper \cite{Per} just cited
develops a-posteriori theorems for state dependent perturbations of periodic orbits 
in ODEs.  In the present work we implement a procedure sufficient for 
verifying the hypotheses of \cite{Per} in concrete examples.
Here we have to balance two competing considerations. 
On the one hand, our arguments require a great deal of quantitative 
information about the perturbing functions \textit{and}
periodic orbit of the unperturbed system.  On the other hand, 
we do not want to restrict our attention to ODEs where explicit formulas for 
periodic orbits are known analytically.  Indeed, our goal is to describe an approach 
which works \textit{in principle} for any ODE with an isolated periodic orbit.  

In the present work these constraints are simultaneously satisfied using computer-assisted
methods of proof for the ODE.  This is a very active area of research, and many viable 
options exist for studying periodic solutions.  A thorough review of the literature is 
a task beyond the scope of this modest introduction, and we refer the interested 
reader to the review articles \cite{MR2652784,VANDENBERG_Dynamics,MR4283203}
and books \cite{TUCKER_ValidatedIntroduction,MR3971222}.
In the present work, for reasons that will be elucidated throughout the manuscript, 
we employ a computer-assisted method of proof wherein  
the ODE is projected into a Banach space of rapidly decaying Chebyshev coefficients.
The truncated problem is solved numerically using Newton's method, 
and the existence of a true solution with Chebyshev
series coefficients near our numerical approximation is proven using a
Newton-Kantorovich argument.  
The approach is adapted from \cite{AJPJ}, 
using techniques also from \cite{JPC,Lessard2018}.  A readable introduction 
to these ideas is found in the first three chapters of \cite{MR3822720}.
Combining the analytical results from \cite{Per} with the mathematically 
rigorous computational methods just discussed, we prove results for a number of  state-dependent delayed perturbations of a van der Pol equation.
As far as we are aware, the present paper and 
\cite{kevinIntegrator} - which deals with the problem of rigorous integration --
are the only two papers thus far in the literature dealing with computer assisted
proofs for SDDEs. 

%
%

%

\begin{remark}[SDDE in electrodynamics] \label{rem:electro}
DDEs arise naturally when modeling distributed 
systems where communication lags between the various subsystems
cannot be ignored.  When the lags themselves depend on the state of the system, we have 
SDDEs.  Examples are common in biology, control theory, and epidemeology, and 
we refer to \cite{MR2457636} for an expansive discussion of relevant 
applications. 
We single out for further discussion an interesting collection of
problems where the perturbative 
arguments developed in the present work could, under suitable modifications, 
be very useful.

While the classical Newtonian theory of 
$N$-body interactions treats gravitational disturbances as acting instantaneously,  it
would be more realistic to incorporate the finite propagation speed of light,
perhaps as a delay.
For example, since the magnitude of the gravitational force due to one body acting on another
is inversely proportional to the square of distance between them, the delay would 
depend on the distance divided by the speed of light.
For systems of particles where this fraction is small, 
this would result in a
state dependent perturbation of the classical equations of motion.

A much more explicit example is found in the 1949 paper of 
Feynmann and Wheeler \cite{MR0032447}, where they put forward a
 theory of direct interparticle action for a system of
point charges in electromagnetic interaction, 
and discuss its potential advantages. 
Subsequent work in this direction is found in Driver's Ph.D. dissertation
on the electromagnetic two body problem \cite{MR2613114}, and in 
several of his subsequent works
\cite{MR0151110,MR0146486,MR0151110,MR240421,MR1065250}.

For a theory based on SDDEs to be symmetric under the reversal of time, it is necessary to 
include an ``advanced'' delay which restores the symmetry of the equations.
This induces a counter intuitive dependence in the equations of motion on the future state of
 the system, in addition to its past
\cite{MR527585,MR553005}.  While dependence on future state is 
an unusual feature of a physical model, there continues to be much interest 
in this approach \cite{MR2639911,MR2912694,MR3389782,MR3456815}.

We mention these works mainly as an 
opportunity to stress that the a-posteriori framework employed in the present work 
allows for the inclusion of advanced delays, and could in principle be applied to the
electromagnetic theories above.  
On the other hand, in their current formulation, the results of \cite{Per} require that 
the periodic orbit of the ODE is isolated, and this never occurs in the Hamiltonian ODEs of 
classical dynamics.  Development of an a-posteriori framework which generalizes the 
approach of \cite{Per} to systems with continuous symmetries is the subject 
of an upcoming work.

\end{remark}


\begin{remark}[Regularity of the results] \label{rem:regularity}
The arguments in the present work are formulated in 
$C^{1+Lip}$ spaces and, though we can bootstrap to obtain more regularity, 
we do not obtain bounds on higher derivatives.  
Bounds on higher derivatives can be obtained using the methods of 
\cite{Per}, but this requires putting more work into the estimates.  

\end{remark}

\begin{remark}[Spectral Bases: Chebyshev versus Taylor] \label{rem:ChebIntro}
The computer-assisted proofs in this paper are formulated 
using a spectral representation --namely Chebyshev series-- for 
the periodic orbit.  For the purposes of this paper it is quite valuable to have a 
representation of the unperturbed periodic solution as an object in $C^k([0,T], \mathbb{R}^d)$ (with 
$k \in \mathbb{N}$, 
or $k = \infty$, or even real analytic), 
where $T$ is the period and $\mathbb{R}^d$ is the state space for the ODE
(rather than for example a representation of the orbit as a fixed or periodic point 
in a Poincar\'e section).
This is because our method requires bounds on the size of the orbit and its derivatives, 
and these bounds are easily recovered from the spectral representation.  

We note that it is also very natural to study periodic solutions using 
Fourier series, and to formulate computer-assisted proofs on Banach spaces 
of rapidly decaying Fourier coefficients as in  \cite{AJPJ}.  However, the 
a-posteriori theory of \cite{Per} requires bounds on solutions of 
some variational equations associated with the periodic orbit, and solutions 
of the variational equation are not periodic.  This problem could be addressed
in Fourier space using the Floquet methods developed  \cite{MR3436565}.  
However, in the present work we found it 
expedient to solve the variational equations directly using a Chebyshev scheme.
This requires a Chebyshev representation of the periodic orbit, which is why we
use Chebyshev series throughout.  
\end{remark}

The paper is organized as follows. In Section~\ref{sec:Formulation}, we recall some 
background from \cite{Per} and introduce a set of polynomial inequalities whose successful 
verification leads to the existence of periodic orbits in the perturbed SDDE. 
In Section~\ref{sec.alg}, we introduce our main algorithm, which provides a way of computing the coefficients of the polynomials and optimizing their variables so that the polynomial inequalities are satisfied.
In Section~\ref{sec.experi}, we apply our algorithm to prove the existence of periodic orbits in a state-dependent delayed perturbation of the Van der Pol equation. We conclude the paper in Section~\ref{sec:conclusion}.

\section{Formulation of the Problem and the Polynomial Inequalities} \label{sec:Formulation}

In this section, we recall some results from \cite{Per} and describe a set of polynomial inequalities whose successful 
verification leads to the existence of the periodic orbits for the perturbed SDDE.  

Consider a smooth ODE on $\R^d$
\begin{equation}\label{ode}
\dot x(t)=f(x(t)),
\end{equation}
and assume it has a periodic orbit $\Gamma$, which is parameterized by $\kzero\colon \T\to \R^d$ where $\T=\R/\Z$. 
In other words, $\Gamma = \{ \kzero(t) : t \in [0,1] \}$.
We consider a perturbation involving a state-dependent (forward or backward) delay term, that is
\begin{equation}\label{sdde}
\dot x(t)=f(x(t))+\ep P\Big(x\big(t-r(x(t))\big)\Big),
\end{equation}
where $r:\R^d\to \R$ is not restricted to be positive. 

When all but one Floquet multipliers of $\Gamma$ are different from $1$, the result from \cite{Per} ensures
that the perturbed equation \eqref{sdde} also has a periodic orbit for
small enough $\ep$.
In order to find values of $\ep$ so that the periodic orbit
persists, we first summarize the proof in \cite{Per}. The proof there
is based on the parameterization method, first introduced in
\cite{MR2177465,MR1976080,MR1976079}.

The parameterization $\kzero\colon \T\to \R^d$ of the periodic orbit for \eqref{ode} with frequency $\om_0$ satisfies
\begin{equation} \label{po}
\omega_0D\kzero(\theta)=f(\kzero(\theta)).
\end{equation}

Consider perturbations $\kh$ and $\omh$, of $\kzero$ and $\om_0$ respectively, so that $K \bydef \kzero+\kh$
parameterizes the periodic orbit of the perturbed equation
\eqref{sdde} and $\om \bydef \om_0+\omh$ is the new frequency. Then $\kh$ and
$\omh$ satisfy
\begin{equation}\label{inveq}
\om_0 D\kh(\theta)-Df(\kzero(\theta))\kh(\theta) =
B^\ep(\theta,\omh,\kh)-\omh D\kzero(\theta),
\end{equation}
where
\begin{align*}
B^\ep(\theta,\omh,\kh)&\bydef N(\theta, \kh)+\ep
P(\kt(\theta))-\omh D\kh(\theta),\\ 
N(\theta, \kh)&\bydef
f(\kzero(\theta)+\kh(\theta))-f(\kzero(\theta))-Df(\kzero(\theta))\kh(\theta), 
\end{align*}
and the term $\kt$ coming from the delay is given by
\[
\kt(\theta)\bydef K(\theta - \om r (K(\theta))). 
\]

To solve equation \eqref{inveq}, we consider the variational
equation. For any fixed $\theta_0\in \T$, let $\Phi(\theta;\theta_0)$
be such that
\begin{equation}\label{var}
\omega_0\frac{d}{d\theta}\Phi(\theta;\theta_0) =
Df(\kzero(\theta))\Phi(\theta;\theta_0),\quad
\Phi(\theta_0;\theta_0)=Id.
\end{equation}
where $Id$ is the identify matrix in $\R^d$.  Then the condition that
$\Gamma$ has all but one Floquet multiplier
different from $1$ is equivalent to the following nondegenerate
assumption \eqref{hypothesisH} on $\Phi(\theta_0+1;\theta_0)$:
\begin{enumerate}
\renewcommand{\theenumi}{H}
\renewcommand{\labelenumi}{(\theenumi)}
\item \label{hypothesisH}\textit{$\Phi(\theta_0+1;\theta_0)$ has a
  simple eigenvalue 1 whose eigenspace is generated by
  $DK_0(\theta_0)$.}
\end{enumerate}

At the point $\kzero(\theta_0)$ on the periodic orbit, the tangent
space has a spectral splitting,
\begin{equation}\label{split}
  T_{\kzero(\theta_0)}\R^d=E_{\theta_0}\oplus Span\{
  D\kzero(\theta_0)\}.
\end{equation} 
We denote the projections onto $Span\{ D\kzero(\theta_0)\}$ and
$E_{\theta_0}$ as $\pt_{\theta_0}$ and $\pp_{\theta_0}$, respectively.

Now, equation \eqref{inveq} can be solved by a fixed point approach
using the variation of constants method. The fixed point of the
following operator $\Gamma^\ep$ solves \eqref{inveq}, where
\[
\Gamma^\ep(\omh,\kh)=\begin{pmatrix}
\Gamma^\ep_1(\omh,\kh)\\ \Gamma^\ep_2(\omh,\kh)
\end{pmatrix},
\]
with
\begin{align*}
 \Gamma^\ep_1(\omh,\kh) & \bydef \frac{\left\langle\int^{\theta_0+1}_{\theta_0}\pt_{\theta_0}\bigl(\Phi(\theta_0+1;s)B^\ep(s,\omh,\kh)\bigr)\, ds,
   D\kzero(\theta_0)\right\rangle}{\left| D\kzero(\theta_0)\right|^2},
\\
\Gamma^\ep_2(\omh,\kh)(\theta) &\bydef \Phi(\theta;\theta_0)u_0
 + \frac{1}{\om_0}\int^{\theta}_{\theta_0}\Phi(\theta;s)\big(B^\ep(s,\omh,\kh)-\Gamma^\ep_1(\omh,\kh)
 D\kzero(s)\big)\, ds
\end{align*}
where $u_0\in E_{\theta_0}$ satisfies 
\[
[Id-\Phi(\theta_0+1;\theta_0)]u_0=\frac{1}{\om_0}\int^{\theta_0+1}_{\theta_0}\pp_{\theta_0}\bigl(\Phi(\theta_0+1;s)B^\ep(s,\kh,\omh)
\bigr) \, ds.
\]
Basically, $\Gamma^\ep_2$ provides the updated $\kh$ with variation of
constants formula. The initial condition $u_0$ is chosen such that
$\Gamma^\ep_2$ is periodic, and $\Gamma^\ep_1$ makes sure that there
exists such $u_0$.

Let $I_a \bydef [-a, a]$, and
\begin{equation}\label{space}
  \mathcal{E}_\beta \bydef \Big\{g\colon \T\to\R^d\big|~ \text{$g$ is
    $C^{2}$}, ~\Big\|\frac{d^i}{d\theta^i}g(\theta)\Big\|\leq
  \beta_i,\ i=0,1,2\Big\},
\end{equation}
where $\| \cdot \|$ is the $C^0$ norm. Note that we could choose any norm in $\R^d$ 
(e.g. Euclidean norm), but a different choice would lead to a slight change in the estimates.
Define the fixed point operator $\Gamma^\ep$ on
$I_a\times \mathcal{E}_\beta$. For small given $\ep$, if one can show
that: (i) $\Gamma^\ep$ maps $I_a\times \mathcal{E}_\beta$ into itself,
and (ii) $\Gamma^\ep$ is a contraction in $C^0$ distance, then there
is a fixed point $(\omh^*,\kh^*)$ of $\Gamma^\ep$. Notice that
$\omh^*\in I_a$ and $\kh^*$ is $C^{1+\Lip}$, that is differentiable with
Lipschitz derivative, using the fact that the closure of $ \mathcal{E}_\beta$
is a subset of $C^{1+\Lip}$ functions. Therefore, equation
\eqref{inveq} is solved. Indeed, (i) ensures the existence of a fixed
point, and (ii) guarantees the uniqueness of the fixed point. As
demonstrated in \cite{Per}, in order to verify (i) and (ii), it is
sufficient to verify the six following inequalities.
\begin{align*}
|\gaep_1(\omh,\kh)|\le a \quad &\Longleftrightarrow \quad Q(\ep, a, \beta_0)\le a\\
\|\gaep_2(\omh,\kh)\|\le \beta_0  \quad& \Longleftrightarrow \quad P_0(\ep, a, \beta_0)\le \beta_0\\
\left\|\frac{d}{d\theta}\gaep_2(\omh,\kh)\right\|\le \beta_1\quad &\Longleftrightarrow \quad P_1(\ep, a, \beta_0, \beta_1)\le \beta_1\\
\left\|\frac{d^2}{d\theta^2}\gaep_2(\omh,\kh)\right\|\le \beta_2 \quad &\Longleftrightarrow \quad P_2(\ep, a, \beta_0, \beta_1,\beta_2)\le \beta_2\\
\gaep \text{ is a contraction in } C^0 \quad &\Longleftrightarrow \quad \begin{cases}
\mu_1(\ep,  a,\beta_0,\beta_1) < 1\\
\mu_2(\ep,  a,\beta_0,\beta_1) < 1
\end{cases}
\end{align*}

As we shall see now in the next section, $Q,~P_0,~P_1,~P_2,~\mu_1$, and $\mu_2$ are polynomials with
computable coefficients which are determined by the unperturbed equation, the
perturbation term $P$, and the delay term $r$.  

\subsection{The explicit construction of the polynomials \boldmath$Q$\unboldmath, \boldmath$P_0$\unboldmath, \boldmath$P_1$\unboldmath, \boldmath$P_2$\unboldmath, \boldmath$\mu_1$\unboldmath~and~\boldmath$\mu_2$\unboldmath} \label{sec.polys}

For a given norm $|\cdot|$ on $\R^d$, denote $|A|$ as the operator
norm of matrix $A$. Define constants
\begin{equation}
\label{C1}
\begin{split}
&C_{1,1} \bydef \|\Phi(\theta_0+1;s)\|=\max_{s\in [\theta_0,\theta_0+1]}\left|\Phi(\theta_0+1;s)\right|,\\
&C_{1,2} \bydef \|\Phi(\theta;s)\|=\max_{\theta\in [\theta_0,\theta_0+1],~s\in [\theta_0,\theta]}\left|\Phi(\theta;s)\right|,\\
&C_{1,3} \bydef \|\Phi(\theta;\theta_0)\|=\max_{\theta\in [\theta_0,\theta_0+1]}\left|\Phi(\theta;\theta_0)\right|.
\end{split}
\end{equation}
Let 
\begin{align*}
&\left\|\frac{d}{ds}\Phi(\theta_0+1;s)\right\|=\max_{s\in [\theta_0,\theta_0+1]}\left|\frac{d}{ds}\Phi(\theta_0+1;s)\right|,\\
&\left\|\frac{d}{ds}\Phi(\theta;s)\right\|=\max_{\theta\in [\theta_0,\theta_0+1],~s\in [\theta_0,\theta]}\left|\frac{d}{ds}\Phi(\theta;s)\right|,
\end{align*}
and
\begin{equation}
\label{C2}
\begin{split}
&C_{2,1} \bydef 1+|\Phi(\theta_0+1;\theta_0)|+\left\|\frac{d}{ds}\Phi(\theta_0+1;s)\right\|,\\
&C_{2,2} \bydef 1+\|\Phi(\theta;\theta_0)\|+\left\|\frac{d}{ds}\Phi(\theta;s)\right\|.
\end{split}
\end{equation}
Recall the spectral splitting in \eqref{split}, let
$\left\|[Id-\Phi(\theta_0+1;\theta_0)]^{-1}\big|_{E_{\theta_0}}\right\|$
be the norm of the operator $[Id-\Phi(\theta_0+1;\theta_0)]^{-1}$
defined in the space $E_{\theta_0}$, let
$\|\pt_{\theta_0}\|=\|\pp_{\theta_0}\|$ be the norm of the projections
$\pt_{\theta_0}$ and $\pp_{\theta_0}$, and let
\begin{equation} \label{eq:constant_M}
M \bydef \frac{\left\|[Id-\Phi(\theta_0+1;\theta_0)]^{-1}\big|_{E_{\theta_0}}\right\|\|\pp_{\theta_0}\|}{\om_0}.
\end{equation}
Now we are ready to provide the polynomials explicitly.
\begin{equation}\label{q}
Q(\ep, a, \beta_0)=\frac{\|\pt_{\theta_0}\|}{\big| D\kzero(\theta_0)\big|}\Big(\ep C_{1,1}\|P\|+C_{2,1}a\beta_0+\frac{C_{1,1}}{2}\|D^2f\|\beta_0^2\Big),
\end{equation}
where $\|P\|$ is the $C^0$ norm of $P$ in a neighborhood of the
periodic orbit of the unperturbed equation \eqref{ode} (containing
the periodic orbit for the perturbed equation \eqref{sdde}), i.e., a
small neighborhood of size $2\beta_0$ around the unperturbed periodic
orbit, and $\|D^2f\|$ is the supremum of the norm of the bilinear
operator $D^2f$ in the aforementioned neighborhood of the unperturbed
periodic orbit, see more details in section \ref{sec:bound_nbhd}. In
the following, all norms without specification mean the supremum norm.
\begin{align}
\nonumber
P_0(\ep, a, \beta_0) & = \ep \left(M C_{1,3} C_{1,1} +\frac{C_{1,2}}{\om_0} \right)\|P\|+\frac{C_{1,2}}{\om_0}\| D\kzero\|a + \left(M C_{1,3} C_{2,1} + \frac{C_{2,2}}{\om_0}\right) a\beta_0 \\
\label{p0}
& \quad  +\frac{\|D^2f\|}{2} \left(M C_{1,3} C_{1,1} +\frac{C_{1,2}}{\om_0} \right)\beta_0^2,
\\
\label{p1}
P_1(\ep, a, \beta_0, \beta_1)&=\frac{1}{\om_0}\Big(\ep\|P\| +\|D\kzero\|a+\|Df\circ\kzero\|\beta_0 + a\beta_1+\frac{1}{2}\|D^2f\|\beta_0^2\Big),
\\
\label{p2}
P_2(\ep, a, \beta_0, \beta_1,\beta_2)& =\frac{1}{\om_0}\Big(\|D^2\kzero\| a+\|D^2f\circ\kzero\|\|D\kzero\|\beta_0 + \|Df\circ\kzero\|\beta_1 +\Big\|\frac{d}{d\theta}B^\ep\Big\|\Big),
\end{align}
where
\begin{multline*}
\Big\|\frac{d}{d\theta}B^\ep\Big\|\le \ep\|DP\|(\| D \kzero \| + \beta_1) [1 + (\om _0 + a) \| Dr \| (\| D \kzero \| + \beta _1)]+a\beta_2\\+\|D^3f\|(\|DK_0\|+\beta_1)\beta_0^2+2\|D^2f\|\beta_0\beta_1.
\end{multline*}
Let
\[
\xi \bydef \ep\|DP\|+\ep\|DP\|(\|D\kzero\|+\beta_1)(\|Dr\|(\om_0+a)+\|r\|)+\|D^2f\|\beta_0.
\]
Then
\begin{equation}\label{mu1}
\mu_1(\ep,  a,\beta_0,\beta_1)=\frac{\|\pt_{\theta_0}\|C_{2,1}}{\big| D\kzero(\theta_0)\big|}(\beta_0+a)+\frac{\|\pt_{\theta_0}\|C_{1,1}}{\big| D\kzero(\theta_0)\big|}\xi, 
\end{equation}
\begin{multline}\label{mu2}
\mu_2(\ep,  a,\beta_0,\beta_1)=\left(M C_{1,3}C_{2,1}+\frac{C_{1,2}\|\pt_{\theta_0}\|C_{2,1}\|D\kzero\|}{\om_0\big| D\kzero(\theta_0)\big|}+\frac{C_{2,2}}{\om_0}\right)(\beta_0+a)\\
+\left( M C_{1,3} C_{1,1}+\frac{C_{1,2}\|\pt_{\theta_0}\|C_{1,1}\|D\kzero\|}{\om_0\big| D\kzero(\theta_0)\big|}+\frac{C_{1,2}}{\om_0}\right)\xi.
\end{multline}
Note that similar to $\|P\|$ and $\|D^2 f\|$, $\|r\|$, $\|Dr\|$, $\|DP\|$, and $\|D^3f\|$ are the supremum
norms in the same neighborhood of the unperturbed periodic orbit mentioned before. 

In the next section, we introduce an algorithm which provides an efficient and automated way of (a) computing the coefficients of the polynomials $Q$, $P_0$, $P_1$, $P_2$, $\mu_1$ and $\mu_2$; and (b) optimizing their variables $\varepsilon$, $a$, $\beta_0$, $\beta_1$ and $\beta_2$ so that the polynomial inequalities are satisfied, hence providing rigorous and constructive proofs of existence of periodic orbits for the perturbed SDDE.

\section{Algorithm} \label{sec.alg}

In the section, we first introduce our main algorithm. We assume that
the $C^1$ norms of the perturbation $P$ and the forward or backward delay
$r$ are given as inputs.

After that, we provide details for each step in the algorithm. In
particular, the computations of the bounds and the constants
appearing in the inequalities are elaborated. In the algorithm, the
word {\em compute} means that we use computer-assisted proofs to
obtain rigorous enclosures of the quantities we are computing.

\newtheorem{alg}[theorem]{Algorithm}
\begin{alg}
\label{alg.delay-cap}
\
\begin{enumerate}
\setlength{\itemsep}{.8em}
\renewcommand*{\theenumi}{\emph{\arabic{enumi}}}
\renewcommand*{\labelenumi}{\theenumi.}
 \item [$\star$] \texttt{Input:} The model $\dot x(t)=f(x(t))$ as in
   \eqref{ode}, $\|P\|$, $\|DP\|$, $\|r\|$, and $\|Dr\|$ for $P$ and
   $r$ in \eqref{sdde}.
 
 \item [$\star$] \texttt{Output:} Positive constants
   $\ep,~a,~\beta_0,~\beta_1,~\beta_2$ such that the inequalities
   \eqref{q}, \eqref{p0}, \eqref{p1}, \eqref{p2}, \eqref{mu1}, and
   \eqref{mu2} are satisfied.

 \item Compute a parameterization of the periodic orbit of
   \eqref{ode}, $\kzero\colon \T\to \R^d$, and the frequency
   $\om_0>0$.

 \item For a fixed $\theta _0 \in \mathbb{T}$, compute the
   solution of the forward variational equation \eqref{var} and also
   the solution of the backward variational equation,
   $\Phi(\theta_0;\theta)=\left(\Phi(\theta;\theta_0)\right)^{-1}$,
   which verifies
   \begin{equation}\label{varb}
     \omega_0\frac{d}{d\theta}\Phi(\theta_0;\theta)=-\Phi(\theta_0;\theta)Df(\kzero(\theta)),
     \quad \Phi(\theta_0;\theta_0)=Id.
   \end{equation}

 \item Compute the eigenvectors of $\Phi(\theta_0
   +1;\theta_0)$.

 \item Compute $\|\pt_{\theta_0}\|$, $\|\pp_{\theta_0}\|$, and $\| [Id
   - \Phi(\theta _0 + 1; \theta _0)]^{-1} | _{E _{\theta _0}} \|$.

 \item Compute the constants $C_{1,1}, C_{1,2}, C_{1,3}$ given in
   \eqref{C1}, $C_{2,1}, C_{2,2}$ given in \eqref{C2}, and $M$ given
   in \eqref{eq:constant_M}.

 \item Optimize $\ep$ for $a,~\beta_0,~\beta_1,~\beta_2$ in certain
   ranges so that the following inequalities are verified.
 \begin{equation}\label{ineqs}
 \begin{cases}
Q(\ep, a, \beta_0)-a\le 0\\
P_0(\ep, a, \beta_0)-\beta_0\le 0\\
P_1(\ep, a, \beta_0, \beta_1)-\beta_1\le 0 \\
P_2(\ep, a, \beta_0, \beta_1,\beta_2)- \beta_2\le 0\\
\mu_1(\ep,  a,\beta_0,\beta_1)-1 < 0\\
\mu_2(\ep,  a,\beta_0,\beta_1) -1< 0.
 \end{cases}
 \end{equation}
 \end{enumerate}
 \end{alg}

\begin{remark}
 Note that Algorithm~\ref{alg.delay-cap} works for a class of
 perturbations $P$ and delays $r$, one only needs the $C^1$ norms of
 them in a neighborhood of the unperturbed periodic orbit.
\end{remark}

\begin{remark}
Note that if $\ep_0$ verifies the inequalities in \eqref{ineqs} given the constants, so will any $0\le\ep\le\ep_0$. 
\end{remark}

\begin{remark}
 There are different ways of optimizing, for example, we can also view
 $\|P\|$, $\| DP \|$, $\| r \|$, and $\| Dr \|$ as variables of the
 polynomials and modify the optimization step in
 Algorithm~\ref{alg.delay-cap} to maximize $\|P\|$, $\| DP \|$, $\| r
 \|$, and $\| Dr \|$ along with $\ep$. The optimization process
 consists, in general, in local searches and thus the initial guesses
 play an important role. The choice of the objective function depends
 on the goals of the problems, which we will specify in
 Section~\ref{sec.experi}.
\end{remark}

\subsection{Computer-Assisted Proofs for the Unperturbed System}
\label{sec.caps-unperturbed}

The first two steps of Algorithm~\ref{alg.delay-cap} require computing
solutions of ODEs, namely a periodic orbit of \eqref{ode} and
solutions of the forward and backward variational equations about the
periodic orbit. These steps are achieved with the tools of rigorously
validated numerics. Using a computer to produce constructive proofs of
existence of solutions of differential equations is by now
well-established, and we refer the interested reader to the survey
papers
\cite{NAKAO_VerifiedPDE,MR2652784,MR1420838,VANDENBERG_Dynamics,GOMEZ_PDESurvey,MR4283203}
and to the books
\cite{MR3971222,TUCKER_ValidatedIntroduction,MR3822720} for more
details. Our approach to compute rigorously the ingredients of Steps 1
and 2 of Algorithm~\ref{alg.delay-cap} uses Chebyshev series expansion
and a Newton-Kantorovich type theorem (i.e. the radii polynomial
approach), as presented in \cite{JPC,AJPJ,MR4292534}. More precisely,
we compute Chebyshev series expansions of the unperturbed periodic
orbit $K _0(\theta)$ of \eqref{po}, the solution $\Phi(\theta;
\theta_0)$ of the forward variational equation \eqref{var}, and the
solution $\Phi(\theta _0;\theta)$ of the backward variational equation
\eqref{varb}. For each of the three problems, a zero-finding problem
is formulated for the Chebyshev coefficients of the solution of the
ODE, which lies in the product of weighed $\ell^1$ spaces that we
denote $\ell_\nu^1$ (for some geometric decay rate $\nu \ge 1$), see
definition in Appendix~\ref{app:l1nu}.

One way to interpret the results is that we have $\bar x$, which is an
approximate finite part of the Chebyshev series, and $\tilde x$, which
is the tail part, such that
\[
 x = \bar x + \tilde x,
\]
although we do not know what is exactly $\tilde x$, we have the bound
\[
\qquad \|\tilde x \|_{\ell ^1_\nu} \leq R,
\] 
for some explicitly given $R>0$, typically quite small.


In order to use Chebyshev series to represent the solutions, it is
standard to rescale the problem and consider solutions defined on the
interval $[-1,1]$, see Appendix~\ref{app.cheb}. Therefore, we define
the scaling parameter $L$, which is related with the period, the
rescaled periodic orbit $O(s)$, the forward variational flow $F(s)$,
and the backward variational flow $B(s)$ as
\begin{equation*}
 \begin{split}
  L &= \tfrac{1}{2\omega _0} , \\
  O(s) &= K _0\bigl(\tfrac{1}{2} (s + 1) \bigr), \\
  F(s) &= \Phi\bigl(\tfrac{1}{2} (s + 1); 0 \bigr), \\
  B(s) &= \Phi\bigl(0; \tfrac{1}{2} (s + 1) \bigr),
 \end{split}
\end{equation*}
for all $s$ in $[-1,1]$, where we note that we fixed $\theta _0 = 0$. In
particular, for all $-1 \leq s, t \leq 1$,
\[
F(t)B(s) =\Phi\bigl(\tfrac{1}{2} (t + 1); 0 \bigr)\Phi\bigl(0;
\tfrac{1}{2} (s + 1) \bigr)=\Phi\bigl(\tfrac{1}{2}(t+1);
\tfrac{1}{2}(s+1)\bigr),
\]
and with this, we can compute $C_{1,2}$ in \eqref{C1} and $C _{2,2}$
in \eqref{C2}.

To solve \eqref{po}, \eqref{var}, and \eqref{varb}, we look for
Chebyshev series expansions of $O(s)$, $F(s)$ and $B(s)$, and each
solution is computed by applying the radii polynomial approach to a
specific zero-finding problem. Assume that this has been achieved,
we have the numerical approximations and estimations of their
tails, that is
\begin{equation}
\label{eq.LOFB}
 \begin{aligned}
  L &= \bar L + \tilde L , & |\tilde L | &\leq R, \\
  O(s) &= \bar O (s) + \tilde O(s), & \| \tilde O \| &\leq R, \\
  F(s) &= \bar F (s) + \tilde F(s), & \| \tilde F \|  &\leq R, \\
  B(s) &= \bar B (s) + \tilde B(s),  & \| \tilde B \| &\leq R,
 \end{aligned}
\end{equation}
in certain norms associated with the norm in $\ell ^1 _\nu$.

\subsection{Computation of the bounds}
\label{sec.compu-bounds}

The bounds in the Algorithm~\ref{alg.delay-cap} require to manage the
information from the Section~\ref{sec.caps-unperturbed} and specify
the norms. From now on, we will stick to Euclidean norm on $\R^d$, then, the supremum norms of vector fields are straightforward.
For the derivatives of the vector fields, we are going to
use the Fr\"obenius norm (i.e. Euclidean norm of the
vectorization) as an upper bound of their operator norms.

Note that $O \colon [-1,1] \to \R ^d$ and $F, B \colon [-1,1] \to \R
^{d \times d}$, and the norms appeared in the coefficients of the
polynomials in Algorithm~\ref{alg.delay-cap} are, in essence, the
supremum norms. Since when $\nu=1$, the $\ell ^1 _\nu$ norm is an upper
bound of the supremum norm, we then have
\begin{equation*}
 \begin{split}
  \| O \| &\leq \sqrt{\sum _{i = 1}^d \| O _i \| _{\ell ^1_1}^2} \leq
  \sqrt{\sum _{i = 1}^d \bigl(\| \bar O _i \| _{\ell ^1_1} + R
   \bigr)^2} , \\ \| F \| &\leq \sqrt{\sum _{i = 1}^d \sum _{j = 1}^d \|
   F_{i,j} \| _{\ell ^1_1}^2} \leq \sqrt{\sum _{i = 1}^d \sum _{j =
      1}^d \bigl(\| \bar F_{i,j} \| _{\ell ^1_1} + R \bigr)^2} ,
 \end{split}
\end{equation*} 
similarly for $B$.

Some of the coefficients in the polynomials involved in
Algorithm~\ref{alg.delay-cap} are now straightforward. However, there
are still a few bounds requiring more computational
effort.

\subsubsection{Bounds on Neighborhoods of the Periodic Orbit}
\label{sec:bound_nbhd}
Some quantities in Algorithm~\ref{alg.delay-cap}, namely $\| D f \|$
appearing in \eqref{p1}, $\| D^2 f\|$ appearing in \eqref{p0},
\eqref{p1} and \eqref{p2}, and $\|D ^3 f\|$ appearing in \eqref{p2},
need to be bounded in a neighborhood of the periodic $K _0$. While
these derivatives are defined everywhere, it is however enough to
consider their bounds in a neighborhood of the unperturbed periodic
orbit.

Therefore, we will rigorously provide a neighborhood enclosing the
unperturbed periodic orbit. Since we use interval arithmetic, the
enclosure will be provided in terms of a hypercube, which admits an
easy computer encoding using interval arithmetic.

Let $O(s)$ be the rigorously proved periodic orbit as in
\eqref{eq.LOFB}. If  $O = (O _1, \dotsc, O _d)$ componentwise, then for each $O_ i$, we consider the optimization
problems of minimization and maximization on $s \in [-1,1]$.

Note that to obtain initial approximations of these optimizations, we
consider a non-interval optimization problem with the numerical
approximations $\bar O _i$, then we verify it using, e.g.,
\texttt{verifyconstraintglobalmin} in \textsc{intlab}, see
\cite{Rump2018}. The rigorous verification provides an interval where
the min/max is located, we then evaluate $\bar O _i$ on this
interval. Now taking into account the errors $\tilde O _i$, we are
ready to provide the infimum and the supremum as boundaries of the
hypercube containing the periodic orbit. Here, one can consider some
safety factors to make the enclosure a little bit bigger
although the process described here already provides a rigorous
enclosure. We then enlarge the hypercube by size $\beta_0$ on the upper and lower bounds, which makes sure that the periodic orbit of the perturbed equation \eqref{sdde} lies in the enlarged neighborhood.

 The outputs of these bounds will be intervals containing the exact bounds,
to prevent a wrapping effect, the hypercube mesh must be adjusted
until those intervals have a small radius. This adjustment will be
model-dependent and often will be an ad-hoc process.

Once the hypercube is determined, we evaluate upper bounds of $\| D f \|$, $\| D^2
f\|$, and $\|D ^3 f\|$ on the hypercube (possibly with a hypercube
mesh) and return the maximum of these evaluations as the bounds. More precisely, for $f = (f _1, \dotsc, f_d)$, using $\partial$ derivative notation, we have that
 on
the hypercube,
\begin{equation*}
  \| D f \| \leq \sqrt{\sum _{i,j = 1} ^d \|\partial _i f _j \|_{\ell^1_1}^2}, \qquad
  \| D ^2 f \| \leq \sqrt{\sum _{i,j,k = 1} ^d \|\partial _{i,j} f _k\|_{\ell^1_1}^2},  \qquad
  \| D ^3 f \| \leq \sqrt{\sum _{i,j,k,l = 1} ^d \|\partial _{i,j,k} f _l\|_{\ell^1_1} ^2},
\end{equation*}
where $\partial _i f _j$ means the partial derivative of $f_j$ with respect to $x_i$, other expressions are similar. Note that for the first inequality above, we used the fact the Fr\"obenius norm is an upper bound of the operator norm of a matrix under Euclidean norm on $\R^d$. The second and third inequalities can be derived from this fact.

\subsubsection{Bounds on the Convolutions} \label{sec:conv}
We have to consider some products involved with the ODE to get several
coefficients for the polynomials in
Algorithm~\ref{alg.delay-cap}. With Chebyshev representations, the
products become convolutions.  Since the $\ell ^1_\nu$ space is a
Banach algebra (see Section~\ref{app:l1nu}), we have $\| a \ast b \|
_{\ell ^1_\nu} \leq \| a \| _{\ell ^1_\nu}\| b\| _{\ell ^1_\nu}$ for
all $a,b \in \ell ^1_\nu$, where $\ast: \ell ^1_\nu \times \ell ^1_\nu
\to \ell ^1_\nu$ denotes the discrete convolution.

However, the above inequality is likely to provide overestimated
bounds which will affect the size of $\varepsilon$ in the optimization
step of the inequalities \eqref{ineqs}. To get better results, we
should avoid using the inequality as much as possible.

More precisely, let $a = \bar a + \tilde a$ and $b = \bar b + \tilde
b$ be $\ell ^1_\nu$ elements with exact truncated parts and the tail
parts. If the tail parts are bounded by $R >0$, then
\begin{equation}
\label{eq.conv-bound}
 \| a \ast b \| _{\ell ^1_\nu} \leq \| \bar a \ast \bar b \| _{\ell
   ^1_\nu} + (\|\bar a \| _{\ell ^1_\nu} + \|\bar b \| _{\ell ^1_\nu}
 + R)R.
\end{equation}

Numerically we keep $\bar a \ast \bar b$ and the bound of the tail
$(\|\bar a \| _{\ell ^1_\nu} + \|\bar b \| _{\ell ^1_\nu} +
R)R$. Thus, we can consider a \texttt{class} that encodes the
truncated Chebyshev series and a bound of its tail. In that class we
overload different operations and make elemental operations, such as
sums, products, norms, easily computable.

Note that the bound in \eqref{eq.conv-bound} becomes more complicated
as we increase the number of the convolutions to bound, i.e. for
cubic, quartic, quintic, etc. convolutions. More precisely, we can
keep the numerical parts and let the tail parts be variables of a
polynomial, e.g. $p(s _1, s _2) = (\bar a + s _1)(\bar b + s _2)$, we
then expand everything in monomials and take into account the bounds
of the tails. If the tails are bounded by $R$, then we let $s _1 = s
_2 = R$ in the expansion.

\subsubsection{Bounds on the Derivatives}
Expressions like $C _{2,1}$ in \eqref{C2} can be bounded by using the
ODE systems. That is,
\[
 \frac{d}{ds} F(1)B(s) = -F(1) L f(O(s)) B(s).
\]
The norm of the righthand side is now easily computable by
convolutions, taking care of the numerical and tail parts of $F(1)$,
$L$, $O(s)$, and $B(s)$.

Another possible way to bound the norm of the derivative of a function
is to use estimates similar to Cauchy bounds.

\subsubsection{Bounds on Triangle Meshes}
The terms $C _{1,2}$ in \eqref{C1} and $ C _{2,2}$ in \eqref{C2} are
computationally expensive because they require considering a
triangular mesh.

Indeed, the terms $F(t) B(s)$ in $C _{1,2}$ and $-F(t) L f(O(s)) B(s)$
in $C _{2,2}$ can be bounded by taking an interval mesh for the
triangle $-1 \leq s \leq t \leq 1$. That is, for a $m(m+1)/2$ mesh
size, we define the intervals
\begin{align*}
 t _k &= -1 + 2 [k-1,k] / m, \\
 s _j &= -1 + 2 [j-1,j] / m,
\end{align*}
for integers $k = 1, \dotsc, m$ and $j = 1, \dotsc, k$. We evaluate
the expressions in these intervals (adding the radius $R$), computing
the norms, and returning the maximum. The value of $m$ is chosen in
such a way that the maximum stagnates with respect to larger $m$.


\subsubsection{Bounds on the Projections}
  
The projections $\Pi _{\theta _0}^\top$ and $\Pi _{\theta _0}^\perp$
have the same norm. To bound them, we first consider the case
$d=2$. The higher dimensional case is similar.

\begin{figure}[ht]
 \begin{center}
  \includegraphics[scale=1]{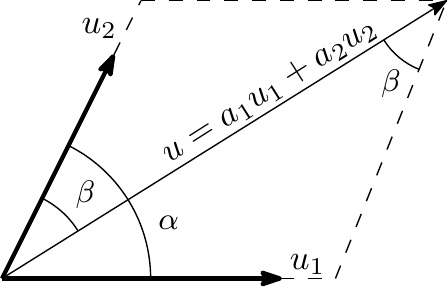}
  \includegraphics[scale=.9]{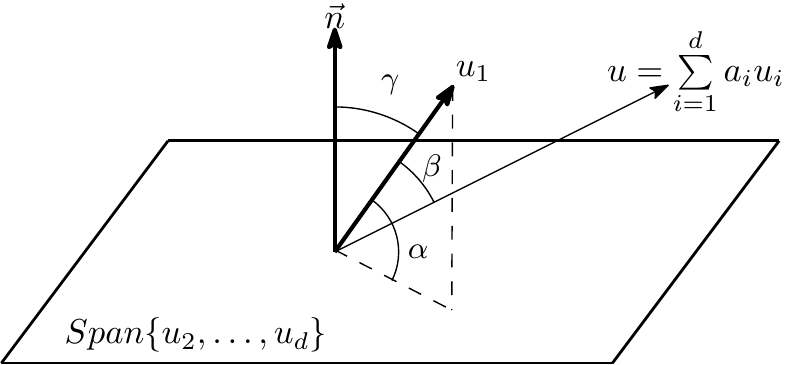}
 \end{center}
 \caption{On the left, the angles in the two-dimensional case $d=2$, a
   vector $u$ in the the span of $u _1$ and $u _2$. On the right, the
   general case, $u _1$ is the eigenvector of the eigenvalue $1$ and
   $\vec n$ is a normal vector of the hyperplane generated by the
   other eigenvectors.}
 \label{fig.angles}
\end{figure}

In the two-dimensional case, the monodromy matrix has two eigenpairs
$(\lambda _i, u _i)$, $i=1,2$, without loss of generality, we assume
that $\lambda _1=1$, $\lambda _2 \ne1$. An arbitrary vector $u$ in the
plane is given by $u = a _1 u _1 + a _2
u _2 $ in the basis $\{u _1, u _2\}$. Let $\alpha$ be the smaller angle between these two
eigenvectors. Then by the trigonometric relations (law of sine) (see Figure~\ref{fig.angles})
\begin{equation*}
 \frac{\| u \|}{\sin \alpha} = \frac{\| a _1 u _1 \|}{\sin \beta} =
 \frac{\| a _2 u _2 \|}{\sin (\alpha - \beta)},
\end{equation*}
then
\[
 \| a _1 u _1 \| = \frac{\sin \beta}{\sin \alpha} \| u \| \qquad
 \text{and} \qquad \| a _2 u _2 \| = \frac{\sin (\alpha - \beta)}{\sin
   \alpha} \| u \|.
\]
As a consequence, the norm of both projections can be bounded as
$\|\Pi _{\theta _0}^\top\|, \|\Pi _{\theta _0}^\perp\| \leq
\frac{1}{\sin \alpha}$.  In practice, the angle $\alpha $ can be
computed by the inner product properties, that is
\begin{equation*}
 \cos \alpha = \frac{u _1 \cdot u _2}{\| u _1 \| \| u _2 \|}.
\end{equation*}
Note that the eigenvectors $u _1$ and $u _2$ are those of the matrix
$F(1)$, which consists of the numerical part and the tail part. These
eigenvectors need to be verified. We can use for example
\texttt{verifyeig} in \textsc{intlab} for the rigorous verification
(see \cite{Rump2001}).

In the $d$-dimensional case, by assumption \eqref{hypothesisH}, the
monodromy matrix $F(1)$ has a simple eigenvalue $1$ with eigenvector
$u_1$. All the other eigenvectors of eigenvalues $\ne 1$ generate a
hyperplane. Let $\alpha$ be the acute angle between $u_1$ and the
hyperplane, then similar to the $d=2$ case, the norms are bounded by
$\frac{1}{\sin \alpha}$. If we consider a normal vector $\vec n$ to
the hyperplane, the angle $\gamma$ between $u_1$ and $\vec n$
satisfies $\sin \alpha=|\cos \gamma|$, see Figure~\ref{fig.angles}.

\subsubsection{Bounds in the complement of tangential directions}

Now we consider
$\left\|[Id-\Phi(\theta_0+1;\theta_0)]^{-1}\big|_{E_{\theta_0}}\right\|$
in Algorithm~\ref{alg.delay-cap}. Due to the hypothesis \eqref{hypothesisH}, we know that restricted to the
subspace $E _{\theta _0}$ of \eqref{split}, the matrix $Id -
\Phi(\theta_0+1;\theta_0)$ is invertible. In our experiments, see
Section~\ref{sec.experi}, since $E_{\theta_0}$ is 1-dimensional, it is easy to invert $Id -
\Phi(\theta_0+1;\theta_0)$ and compute its norm in $E_{\theta_0}$. In general, one could consider the nonzero singular values of $Id -\Phi(\theta_0+1;\theta_0)$. Since under Euclidean norm, the matrix operator norm is the largest singular value, the norm $\left\|[Id-\Phi(\theta_0+1;\theta_0)]^{-1}\big|_{E_{\theta_0}}\right\|$ we want is the reciprocal of the smallest modulus of the nonzero singular values of $Id -\Phi(\theta_0+1;\theta_0)$.

\subsection{Solving the Inequalities}
The last step in Algorithm~\ref{alg.delay-cap} consists in optimizing
the inequalities \eqref{ineqs} such that, for instance, the
perturbative parameter $\varepsilon$ is maximized. The variables of
that optimization are $\varepsilon$, $a$, $\beta _0$, $\beta _1$, and
$\beta _2$. All of them must be strictly positive, that is more than
the epsilon machine, e.g. $2^{-52}$.  Moreover, the result in
\cite{Per} says that $a$ and $\beta _0$ are of the same order as
$\varepsilon$ asymptotically.

The optimization problem that maximizes the $\varepsilon$ is computed
numerically. That is, we take the upper bounds of the coefficient
intervals of the polynomials in \eqref{ineqs}, we apply the numerical
maximization problem, and then we check the result with the
inequalities having the interval coefficients.

A similar process is applied when we try to optimize the class
of perturbations and forward/backward delays, namely $\| P \|$, $\| DP \|$, $\|r \|$,
and $\|Dr\|$. In this situation, one needs to play more
with the objective function of interest, which depends on the goals of
the proof itself.

To prevent a loss of information, it is convenient to expand the
polynomial in \eqref{ineqs} in monomials. Thus the coefficients in the
numerical search is as sharp as possible. Appendix~\ref{app:monopoly}
shows the explicit monomial expressions for the two
different optimizations proposed here.


In the experiments we use the Optimization Toolbox in Matlab 
to get an approximation of the maximizations and \textsc{intlab} to
certify that the validity of the inequalities with the intervals.

\section{The Van der Pol Example}
\label{sec.experi}

For the unperturbed ODE, we consider the Van der Pol equation with
parameter $\mu>0$
\begin{equation}\label{vdp}
 \begin{split}
  \dot x_1 &= x_2 \\
  \dot x_2 &= \mu (1-x_1^2)x_2-x_1.
 \end{split}
\end{equation}

Recall that Algorithm \ref{alg.delay-cap} only requires the $C^1$
bounds of $P$ and $r$ on a bounded set. Hence, we will not fix $P$ and
$r$, but rather work with a class of $P$'s and $r$'s.

\subsection{Results}
The main results consists of three theorems, namely
Theorem~\ref{thm.vdp-ars}, Theorem~\ref{thm.vdp-sdde0} and
Theorem~\ref{thm.vdp-sdde}. The first theorem computes and proves the
existence of the periodic orbit, the forward, and the backward
variational flows of the ODE in \eqref{vdp}. It uses the radii
polynomial approach, see Appendix~\ref{app.radii}, which provides and
explicit distance between the numerical approximation and the exact
solution. The second theorem proves the existence of periodic orbits
in state-dependent delay perturbations of \eqref{vdp} (in the form of
\eqref{sdde}) for some values of the perturbation parameter $\ep$, given the norms of the perturbations and the forward or
backward delays. Moreover, we have estimations on the differences
between the frequencies and the periodic orbits before and after the
perturbation. The third theorem establishes that the perturbation and delay terms can be more general.

\begin{theorem} \label{thm.vdp-ars}
 Fix a parameter value $\mu \in \{ j/10:j=1,\dots,10\}$.  Let $\bar O
 \colon [-1,1] \to \R ^2$ be a numerical solution of the periodic
 orbit of \eqref{vdp} (passing close to the point $(2,0) \in \R^2$)
 consisting of $200$ Chebyshev coefficients (coordinatewise), and
 $2\bar L$ be the numerical approximate period. Then the true periodic
 orbit $O \colon [-1,1] \to \R ^2$ and its period $2L$ satisfy
 \[
  \| O - \bar O \| 
  \leq r _0, \qquad 
  | L - \bar L | \leq r _0
 \]
 where $r _0$ depends on $\mu$ and is given in
 Table~\ref{tab.vdp-ars}.  Moreover, let $\bar F, \bar B \colon [-1,1]
 \to \R ^{2\times 2}$ be numerical approximations of the forward and
 backward variational flows of $O$ represented by $200$ Chebyshev
 coefficients (entrywise). Then the true solutions $F, B \colon [-1,1]
 \to \R ^{2\times 2}$ satisfy
 \[
  \| F - \bar F \| 
  \leq r _1, \qquad 
  \| B - \bar B \|
  \leq r _2
 \]
 where $r _1$ and $r _2$ are given in Table~\ref{tab.vdp-ars}.
 \begin{table}[ht]
  \centering\tiny {\tt
  \begin{tabular}{c||*3c}
   $\mu$ & $ r _0$ & $ r _1$ & $ r _2$ \\ \hline \hline 
0.1 &  5.573887884260317e-13 &  2.642512865973085e-12 &  2.996819243277690e-12 \\
0.2 &  8.709167228557321e-13 &  5.663812970926017e-12 &  1.015950775601682e-11 \\
0.3 &  1.189397347751034e-12 &  1.077390609220761e-11 &  3.082429090530088e-11 \\
0.4 &  1.557869751645740e-12 &  1.959454762049424e-11 &  9.753694436554251e-11 \\
0.5 &  1.928926778070392e-12 &  3.375851743576695e-11 &  3.280028209406751e-10 \\
0.6 &  2.366661800432010e-12 &  5.761835216827375e-11 &  1.155963389952129e-09 \\
0.7 &  2.865419317705822e-12 &  9.665742902893244e-11 &  4.426491597030142e-09 \\
0.8 &  3.544996952415318e-12 &  1.634487711041929e-10 &  1.894089923170409e-08 \\
0.9 &  4.539421525763888e-12 &  2.781116424781816e-10 &  9.518421901009437e-08 \\
1.0 &  5.835732296028395e-12 &  4.555755576317590e-10 &  5.724940629447956e-07 
  \end{tabular}
  }
  \caption{Radii for the computer-assisted proofs for the periodic
    orbit and its period ($r_0$), the forward ($r_1$), and the
    backward ($r_2$) variational flows for different parameter values
    $\mu$ of \eqref{vdp}.}\label{tab.vdp-ars}
 \end{table}
\end{theorem}

\begin{remark}
 The norm of the difference between $O$ and $\bar O$ in above theorem
 should be interpreted as the maximum of the componentwise $\ell
 ^1_{\nu }$ norm for $\nu=1.01$, where each component of $\bar O$
 ($200$ Chebyshev coefficients) is viewed as an element in $\ell
 ^1_{\nu }$ with zero tail. Similar for $ \| F - \bar F \|$ and $\| B
 - \bar B \|$.
\end{remark}

\begin{remark}
 The number $200$ of Chebyshev coefficients in
 Theorem~\ref{thm.vdp-ars} are determined by plotting the coefficients
 and truncating before the stagnation in the tail (due to the
 precision of the arithmetic as $2^{-52}$ in double-precision.)
\end{remark}

Now we are ready to provide our first persistence result.

\begin{theorem}  \label{thm.vdp-sdde0}
 Consider a perturbation of the form \eqref{sdde} to equation
 \eqref{vdp}, assume that the functions $P$ and $r$ satisfy $\|P\|$,
 $\|DP\|$, $\|r\|$, and $\|Dr\|\le 1$, then for values of the
 parameter $\mu$ and constants $a$, $ \beta _0$, $ \beta _1$, and $
 \beta _2$ in the Table~\ref{tab.abep} below, if $ \ep\le\ep _0$ in
 the Table~\ref{tab.abep}, the inequalities in \eqref{ineqs} are
 satisfied. Hence, there exists a $C^{1+Lip}$ periodic orbit for the
 perturbed system \eqref{sdde}. The differences of the unperturbed and
 perturbed periodic orbits and their frequencies lie in the $C^0$
 closure of the space $I _a \times \mathcal{E}_\beta$ as in
 \eqref{space}.
\begin{table}[ht]
\centering\tiny {\tt
\begin{tabular}{c||*5c}
$\mu$ & $a$ & $ \beta _0$ & $ \beta _1$ & $ \beta _2$ & $ \ep _0$ \\ \hline \hline 
0.1 &  3.9289483314e-03 &  1.3814767259e-02 &  1.9351897706e-02 &  3.3489862240 &  2.0479538526e-03 \\
    0.2 &  3.9458701226e-03 &  1.2482231722e-02 &  2.4285482147e-02 &  3.3566339504 &  1.8547822608e-03 \\
    0.3 &  2.9095469512e-03 &  9.9870189178e-03 &  2.8450975266e-02 &  3.3611368094 &  1.1932985497e-03 \\
    0.4 &  1.8632478932e-03 &  7.4544551083e-03 &  4.5960969418e-02 &  3.3707741174 &  6.4266953301e-04 \\
    0.5 &  1.0872644592e-03 &  5.2329709138e-03 &  1.0817095356e-01 &  3.4113619375 &  3.1004411778e-04 \\
    0.6 &  5.9335608942e-04 &  3.5047699947e-03 &  2.8241308709e-01 &  3.5732825217 &  1.3957190885e-04 \\
    0.7 &  3.0480495862e-04 &  2.2543350887e-03 &  4.5110477848e-01 &  3.8463306657 &  5.8261321498e-05 \\
    0.8 &  4.7155981698e-05 &  7.9828837072e-04 &  4.6510209568e-01 &  3.9780706528 &  6.9961093824e-06 \\
    0.9 &  3.2308284738e-06 &  1.8561771388e-04 &  4.0830621482e-01 &  3.9828277487 &  2.3650493428e-07 \\
    1.0 &  2.0172795949e-07 &  4.0905441781e-05 &  3.5946080062e-01 &  3.9846796546 &  9.7963190835e-09
\end{tabular}
}
  \caption{Values of $a$ and $\beta$'s, and admissible $\ep_0$ so
    that the periodic orbit persists under perturbation.}
  \label{tab.abep}
 \end{table}
\end{theorem}

\begin{remark}
 Theorem~\ref{thm.vdp-sdde0} and bootstrapping techniques lead to
 infinite regularity provided that the unperturbed system and the
 functions $P$ and $r$ are all smooth.
\end{remark}

The other parameters admits several options of optimization depending
on the aims.  Notice that the norm $\| P \|$ can be normalized to $1$
since $P$ is multiplied by the perturbative parameter
$\varepsilon$. Our second perturbative result,
Theorem~\ref{thm.vdp-sdde}, optimizes the threshold $\varepsilon _0$
and also the class of perturbation $\| DP \|$ and of the delay $\| r
\|$, $\| Dr \|$. Then for all $\varepsilon \leq \varepsilon _0$, $\| r
\| \leq c _r$, $\|D r\| \leq c _{Dr}$, and $\| DP \| \leq c _{DP}$
with $\varepsilon_0$, $c _r$, $c _{Dr}$, and $c_{DP}$ in
Table~\ref{tab.vdp-values}, the inequalities in \eqref{ineqs} are
satisfied, and so there is a $C^{1+Lip}$ periodic orbit of the
perturbed system.

\begin{theorem} \label{thm.vdp-sdde}
 Let \eqref{vdp} be the unpertubed ODE and consider a perturbation of
 the form \eqref{sdde}. Given $a$, $\beta _0$, $\beta _1$, and $\beta
 _2$ in Table~\ref{tab.vdp-abetas}. Then for all $\varepsilon \leq
 \varepsilon _0$, $\| r \| \leq c _r$, $\|D r\| \leq c _{Dr}$, and $\|
 DP \| \leq c _{DP}$ with $\varepsilon _0$, $c _r$, $c _{Dr}$, and $c
 _{DP}$ in Table~\ref{tab.vdp-values}, the inequalities in
 \eqref{ineqs} are satisfied, and so there exits a $C^{1+Lip}$
 periodic orbit for the perturbed system.

 \begin{table}[ht]
  \centering\tiny{ \tt
  \begin{tabular}{c||*4c}
   $\mu$ & $ a$ & $ \beta _0$ & $ \beta _1$ & $ \beta _2$ \\ \hline \hline
0.1 &  5.3516317207e-03 &  1.3512619469e-02 &  1.8956002692e+00 &  3.7809698565 \\
0.2 &  4.4313887934e-03 &  1.2622565762e-02 &  1.4920065531e+00 &  3.8602413976 \\
0.3 &  2.8488130419e-03 &  9.6105596050e-03 &  1.1632362043e+00 &  3.9047783796 \\
0.4 &  1.6254775156e-03 &  6.6908933879e-03 &  9.2220210133e-01 &  3.9284935852 \\
0.5 &  8.6242391376e-04 &  4.4107683210e-03 &  7.5126159361e-01 &  3.9475219696 \\
0.6 &  4.1234067787e-04 &  2.6447824466e-03 &  6.3209724862e-01 &  3.9667717362 \\
0.7 &  1.9898994536e-04 &  1.6282965579e-03 &  5.3591862272e-01 &  3.9728038918 \\
0.8 &  2.8152415098e-05 &  5.8068726773e-04 &  4.6148172152e-01 &  3.9758037076 \\
0.9 &  2.3423038382e-06 &  1.5112884117e-04 &  4.0304855658e-01 &  3.9667079412 \\
1.0 &  1.5346224834e-07 &  3.3982751869e-05 &  3.6410084636e-01 &  4.0500566273
  \end{tabular}
  }
  \caption{Values of the space $I _a \times \mathcal{E}_\beta$ in
    \eqref{space} for which the perturbed periodic orbit
    exists.}\label{tab.vdp-abetas}
 \end{table}
 
 \begin{table}[ht]
  \centering\tiny{ \tt
  \begin{tabular}{c||*4c}
   $\mu$ & $ \varepsilon _0$ & $ c _r $ & $ c _{Dr}$ & $c_{DP}$ \\ \hline \hline
0.1 &  5.2514914300e-07 &  1.0525530335e+02 &  9.7470170038e+01 &  1.1923883191e+01 \\
0.2 &  4.8780187758e-07 &  1.0550841730e+02 &  9.7829645680e+01 &  1.2607166526e+01 \\
0.3 &  3.5942433417e-07 &  1.0677903035e+02 &  9.7626077504e+01 &  1.3098287211e+01 \\
0.4 &  2.0629988060e-07 &  1.1537382755e+02 &  9.7823817937e+01 &  1.4572329081e+01 \\
0.5 &  9.4989481038e-08 &  1.3431486893e+02 &  7.9752272527e+01 &  2.1261206986e+01 \\
0.6 &  1.0511326205e-07 &  1.5922648281e+02 &  5.4051960357e+01 &  2.1436738121e+01 \\
0.7 &  9.9701256277e-08 &  1.4962219117e+02 &  5.4940586398e+01 &  1.3762820031e+01 \\
0.8 &  3.4717672547e-08 &  2.8352283832e+02 &  1.2981903915e+02 &  5.0876339676e+00 \\
0.9 &  9.0478966022e-09 &  3.9078461001e+02 &  1.1719345674e+02 &  4.8009874216e+00 \\
1.0 &  1.8649496078e-09 &  3.6114164131e+02 &  1.2261159547e+02 &  4.7222224864e+00 
  \end{tabular}
  }
  \caption{Optimized perturbative parameters $\varepsilon _0$, $c _r$,
    $c _{Dr}$, and $c _{DP}$.} \label{tab.vdp-values}
 \end{table}
\end{theorem}

\begin{remark}
 The quantities in Table~\ref{tab.vdp-values} have been obtained by an
 optimization process of the inequalities \eqref{ineqs} with the
 objective function $\ep ^2 \|r \| \|Dr\| \| DP \|^2$ for each of the
 $\mu$ values.
\end{remark}

\subsection{The Zero Finding-problems for the Unperturbed System}
In order to calculate the required coefficients for this example
rigorously, we follow the radii polynomial approach, see
\cite{AJPJ,JPC} and a summary in Appendix \ref{app.radii}, and we use
the \textsc{intlab} package in Matlab, \cite{Ru99a}. We solve a
boundary value problem to get the periodic orbit, and initial value
problems to get the solutions of the forward and backward variational
equations. In the following subsections, each of these problems will
be formulated as a zero-finding problem.

In Algorithm~\ref{alg.delay-cap}, we need to solve for the periodic
orbit $O(s) = (a _1(s), a _2(s))$ and the scaling parameter $L$ that
verifies

\begin{equation}
\label{vdpper}
 \frac{d}{d s}
 \begin{pmatrix}
  a _1 \\ a _2
 \end{pmatrix} = L 
 \begin{pmatrix}
  a _2 \\ \mu (1-a_1^2)a_2-a_1
 \end{pmatrix} , \quad
O(-1)=O(1)  ,\quad\text{and} \quad a _1(-1) = 0,
\end{equation}
the forward system $F(s) = 
\left(
\begin{smallmatrix}
 v _{11} & v _{12} \\ v _{21} & v _{22}
\end{smallmatrix}
\right)(s)$ that verifies
\begin{equation}
\label{vdpfwd}
 \frac{d}{ds}
 \begin{pmatrix}
  v _{11} & v _{12} \\ v _{21} & v _{22}
 \end{pmatrix} = 
 L\begin{pmatrix}
  0 & 1 \\ \vartheta & \varrho
 \end{pmatrix}
 \begin{pmatrix}
  v _{11} & v _{12} \\ v _{21} & v _{22}
 \end{pmatrix} \quad \text{and} \quad F(-1) = Id_2,
\end{equation}
with $\vartheta(s) = -2\mu a _1 a _2 - 1 $ and $\varrho(s) = \mu(1 -
a_1^2)$, and the backward system $B(s) = \left(
\begin{smallmatrix}
 u _{11} & u _{12} \\ u _{21} & u _{22}
\end{smallmatrix}
\right)(s)$ that verifies
\begin{equation*}
 \frac{d}{ds}
 \begin{pmatrix}
 u _{11} & u _{12} \\ u _{21} & u _{22}
 \end{pmatrix} = 
 -L\begin{pmatrix}
 u _{11} & u _{12} \\ u _{21} & u _{22}
 \end{pmatrix}
 \begin{pmatrix}
  0 & 1 \\  \vartheta & \varrho
 \end{pmatrix} \quad \text{and} \quad B(-1) = Id_2.
\end{equation*} 

\begin{figure}[ht]
 \centering
 \includegraphics[scale=.41]{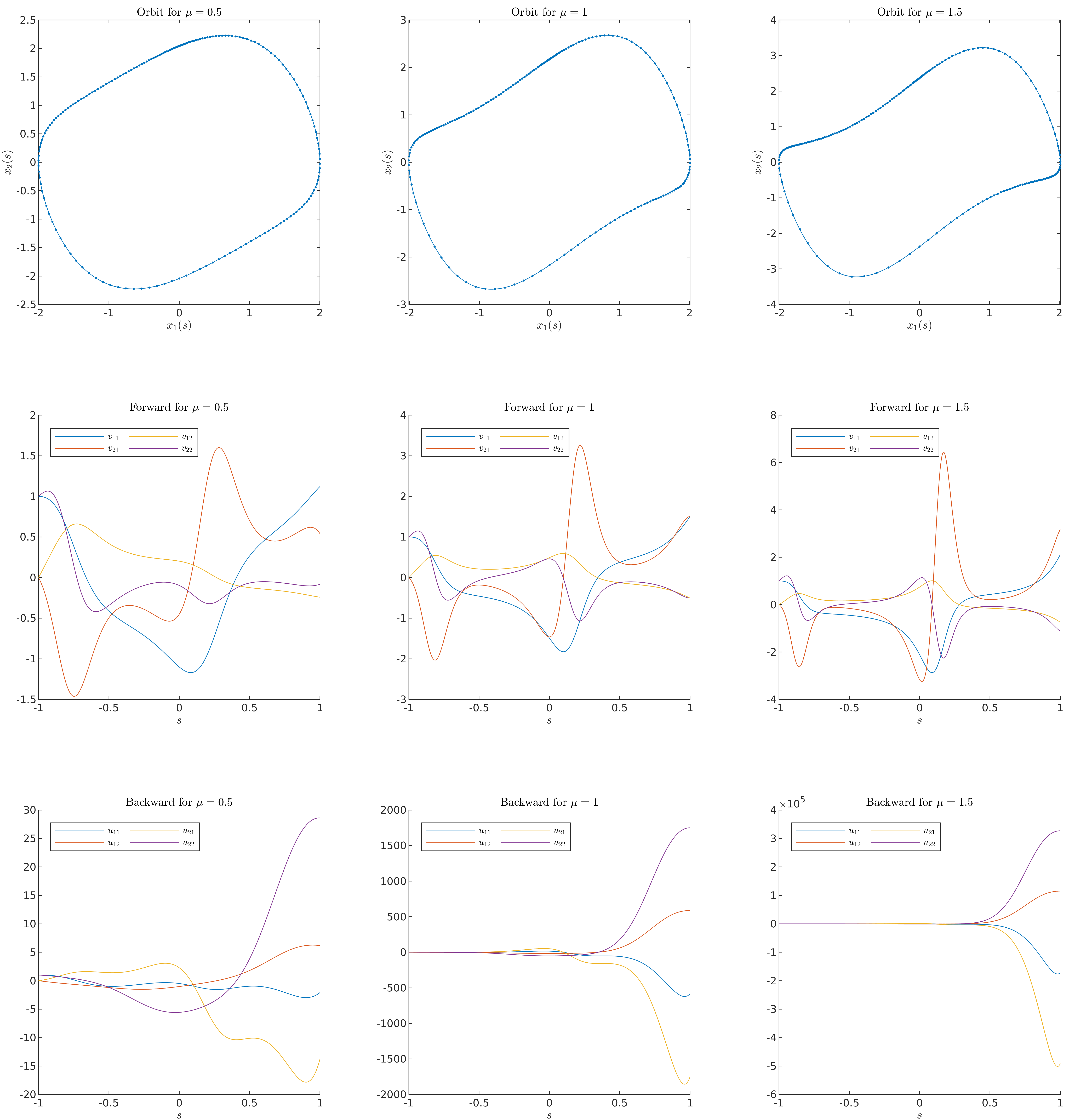}
 \caption{On the left $\mu=0.5$, on the middle $\mu =1$, and on the
   left $\mu = 1.5$ with $\mu$ parameter in \eqref{vdp}.  From top to
   bottom the orbit of the periodic orbit passing close to the point
   $(2,0)$, its forward and backward variational flows.}
 \label{fig.flows}
\end{figure}

We utilize Chebyshev discretization to represent the periodic orbit
and its variational (forward and backward)
flows. Figure~\ref{fig.flows} shows the (numerical) solutions for some
values of the parameter $\mu$. Note that for $\mu \geq 1$ the backward
flow grows rapidly since the periodic orbit attracts strongly as we
can see in its (numerical) real non-trivial eigenvalue $\lambda _\mu$
of the monodromy matrix
\begin{equation*}
 \begin{split}
  \lambda _{0.5} &= \mathtt{3.917692025927352 \times 10^{-2}}, \\
  \lambda _{1}   &= \mathtt{8.596950636046152\times 10^{-4}}, \\
  \lambda _{1.5} &= \mathtt{6.466756568013210\times 10^{-6}}.
 \end{split}
\end{equation*}

\begin{remark}
 Note that in \eqref{vdpper} we considered the spatial section $O
 _1(-1) = 0$. Thanks to hypothesis \eqref{hypothesisH}, we can solve
 the problem successfully.
\end{remark}

\subsubsection{Zero-finding problem for the periodic orbit}
The rigorous proof concerning $L >0$ and the periodic orbit $O \colon
[-1,1] \to \mathbb{R}^2$, with sequence of Chebyshev coefficients
$(a_1, a_2)$, is done in the Banach space $X _\nu \bydef (\ell
^1_\nu)^2 \times \R$. The space is endowed with the product norm
\begin{equation*}
    \| (a_1,a_2,L) \| \bydef \max \{ \| a_1 \| _{\ell ^1_\nu}, \|
    a_2 \| _{\ell ^1_\nu}, |L | \},
\end{equation*}
for $(a_1,a_2,L) \in X _\nu$.

The proof is obtained by solving a boundary value problem, see
Appendix~\ref{app.bvp}, with boundary conditions $a _1(-1) = 0$ (the
phase condition) and $a_i(-1)=a_i(1)$, for $i=1,2$. To define the
equivalent zero-finding map $\mathcal O \colon X _\nu \to X _{\nu'}$,
with $1 < \nu ' \leq \nu $, of \eqref{vdpper}, we first consider
linear operators $M, T \colon \ell ^1_\nu \to \ell ^1_\nu$, and
$\Lambda \colon \ell ^1_\nu \to \ell ^1_{\nu'}$ defined as
\begin{align}
\nonumber
 M &\bydef \begin{pmatrix}
 0 & 1 & 0 & 1 & 0& 1 & 0 & \cdots \\
 0 & 0 & 0 & 0 & 0 & 0 & \cdots \\
 & \ddots & \ddots &  \ddots & \ddots & \ddots \\
 & \cdots & 0 & 0 & 0 & 0 & 0  & \cdots \\
 & & \cdots &  \ddots & \ddots & \ddots & \ddots 
\end{pmatrix}, \\
    \label{eq.T}
T &\bydef \begin{pmatrix}
 0&  0 & 0 & 0 & 0 & 0 & \cdots \\
 -1 & 0 & 1 & 0 & 0 & \cdots \\
 0 & -1 & 0 & 1 &0 & 0 & \cdots \\
 & \ddots & \ddots &  \ddots & \ddots & \ddots \\
 & \cdots & 0 & -1 & 0 & 1 & \cdots \\
 & & \cdots &  \ddots & \ddots & \ddots & \ddots
\end{pmatrix},
\intertext{and }
    \label{eq.Lambda}
 \Lambda &\bydef \begin{pmatrix}
 0 & 0 & 0 & 0 & 0 & 0 & \cdots \\
 0 & 2 & 0 & 0 & 0 & \cdots \\
 0 & 0 & 4 & 0 & 0 & 0 & \cdots \\
 & \ddots & \ddots &  \ddots & \ddots & \ddots & \cdots \\
 & \cdots & 0 & 0 & 2 k  & 0 & \cdots \\
 & & \cdots &  \ddots & \ddots & \ddots
\end{pmatrix}.
\end{align}
Then the zero-finding map for the periodic orbit is given by $\mathcal
O=(\mathcal O _1, \mathcal O _2, \mathcal O _3)$, where
\begin{equation*}
 \begin{aligned}
    \mathcal O _i(a_1,a_2,L) & \bydef (M+\Lambda) a _i +L\cdot T
    f_i(a), & i &\in \{1,2\} ,\\ \mathcal O _3(a_1,a_2,L) & \bydef
    (a_1) _0 + 2\sum _{k \geq 1} (a _1)_k (-1)^k,
 \end{aligned}
\end{equation*}
with
\begin{equation*}
  \begin{split}
    f_1(a) &\bydef   a _2, \\
    f_2(a) &\bydef \mu a _2 - \mu (a _1 * a_1 * a _2) - a _1,
  \end{split}
\end{equation*}
with $\ast$ being the convolution in $\ell ^1_\nu$, see
Appendix~\ref{app.radii}.

Following the radii polynomial approach, we first compute a numerical
solution $(\bar a_1, \bar a_2,\bar L) \in X_\nu$ such that $\mathcal
O(\bar a_1,\bar a_2,\bar L)\approx0$, that is it is zero up to a given
tolerance, for instance $10^{-12}$.

Then computing the radii polynomial \eqref{eq:radii_polynomial} as in
Lemma \ref{lem:radii}, there is an exact solution $(a_1,a_2,L) \in
X_\nu$ within distance $r _0 > 0$ of the numerical solution such that
$\mathcal O(a_1,a_2,L)=0$.

\subsubsection{Zero-finding problem for the solution of the forward variational equation} 
\label{sec:forward}

To compute the solution of \eqref{vdpfwd}, we first assume that one
computed rigorously the values of $L$ and the periodic orbit $ (a_1,
a_2)$ verifying \eqref{vdpper}. We split these information as the sum
of the numerical parts and the tail parts, that is
\begin{equation*}
a_i=\bar a_i+\tilde a_i, \quad \text{for }i \in \{1,2\}, \qquad 
L=\bar L+\tilde L.
\end{equation*}
From the proof, we know that $\|\tilde a_i\| _{\ell ^1_\nu}\le r_0$
and $|\tilde L| \le r_0$.

The $\vartheta$ and $\varrho$ in \eqref{vdpfwd} depend on the periodic
orbit, and are represented as Chebyshev series. We split them by their
numerical and tail parts
\[
\vartheta \bydef \bar \vartheta + \tilde \vartheta, \qquad 
  \varrho \bydef \bar \varrho + \tilde \varrho.
\]
Explicitly, the numerical parts are given by
\begin{equation*}
 \bar \vartheta _k \bydef -2 \mu (\bar a _1 \ast \bar a _2)_k -
 \delta _{0,k}, \qquad \bar \varrho _k \bydef \mu \delta _{0,k} -
 \mu (\bar a _1 \ast \bar a _1)_k,
\end{equation*}
and the tails contain the crossing terms from the convolutions.

Let $S \colon \ell^1_\nu \to \ell^1_\nu$ be the linear operator
defined by
\begin{equation*}
    S \bydef
    \begin{pmatrix}
    -1 & 2 & -2 & 2 & -2 & \cdots \\
    0  & 0 &  0 & 0 &  0 & \cdots \\
    \vdots & \vdots & \vdots & \vdots & \vdots & \ddots \\
    0  & 0 &  0 & 0 &  0 & \cdots \\
    \end{pmatrix}
\end{equation*}
and let $\ell _\alpha \colon \ell ^1_\nu \to \ell ^1_\nu$ be the
linear operator of convolution multiplication defined, given $\alpha
\in \ell ^1_\nu$, by $\ell _\alpha \colon b \mapsto \alpha \ast b$.

The zero-finding problem equivalent to \eqref{vdpfwd} consists in two
proofs of initial value problems, see Appendix~\ref{app.ivp}. Each of
them corresponds to a column of $F(s)$, and they are given by the map
$\mathcal F _p \colon (\ell ^1_\nu) ^2 \to (\ell ^1_{\nu'})^2$, with
$1 \leq \nu ' < \nu$ defined as
\begin{equation}
\label{eq.zeroopfwd}
    \begin{split}
        \mathcal F _p (v) &\bydef \bar {\mathcal F} _p (v) + \tilde
                 {\mathcal F} _p(v), \\ \bar {\mathcal F} _p(v)
                 &\bydef p +
        \begin{pmatrix}
           S  + \Lambda & 0 \\
           0 & S  + \Lambda
        \end{pmatrix} v + \bar L  \begin{pmatrix}
            0 & T \\ T \ell _{\bar \vartheta} & T \ell _{\bar \varrho}
        \end{pmatrix} v, \\
        \tilde {\mathcal F} _p(v) &\bydef  \tilde L \begin{pmatrix}
            0 & T \\ T \ell _{\bar \vartheta} & T \ell _{\bar \varrho}
        \end{pmatrix} v + (\bar L + \tilde L)  \begin{pmatrix}
            0 & 0 \\ T\ell _{\tilde \vartheta} & T\ell _{\tilde \varrho}
        \end{pmatrix} v,
    \end{split}
\end{equation}
where the operators $T$ and $\Lambda$ are defined in \eqref{eq.T} and
\eqref{eq.Lambda} respectively.

The element $p\in (\ell^1_\nu)^2$ in \eqref{eq.zeroopfwd} is related
to the initial condition of the chosen column, when $v$ corresponds to
the first column of $F(s)$, then $p$ is the vector of Chebyshev
coefficients of the constant vector
$\left(\begin{smallmatrix}1\\0\end{smallmatrix}\right)$; otherwise,
  $p$ is the vector of Chebyshev coefficients of the constant vector
  $\left(\begin{smallmatrix}0\\1\end{smallmatrix}\right)$.
  
Now computing the elements in Lemma~\ref{lem:radii} for the initial
conditions, we prove the existence of the solution $F(s) = \bar F(s) +
\tilde F(s)$ with $\| \tilde F \| = \max _{i,j \in \{1,2\}} \| \tilde
v _{ij}\| _{\ell ^1_\nu} \leq r _1$.

\subsubsection{Zero-finding problem for the solution of the backward variational equation}
Similarly to the forward flow in \eqref{eq.zeroopfwd}, we consider the
map $\mathcal B _p \colon (\ell ^1_\nu) ^2 \to (\ell ^1_{\nu'})^2$
defined by
\begin{equation}
\label{eq.zeroopbwd}
    \begin{split}
        \mathcal B _p (u) &\bydef  \bar {\mathcal B} _p (u) + \tilde {\mathcal B} _p(u),  \\
        \bar {\mathcal B} _p(u) &\bydef p + 
        \begin{pmatrix}
           S  + \Lambda & 0 \\
           0 & S  + \Lambda
        \end{pmatrix} u - \bar L  \begin{pmatrix}
            0 & T \ell _{\bar \vartheta} \\ T & T \ell _{\bar \varrho}
        \end{pmatrix} u, \\
        \tilde {\mathcal B} _p(u) &\bydef - \tilde L \begin{pmatrix}
            0 & T \ell _{\bar \vartheta} \\ T & T \ell _{\bar \varrho}
        \end{pmatrix} u - (\bar L + \tilde L) \begin{pmatrix}
            0 & T \ell _{\tilde \vartheta}  \\ 0 & T\ell _{\tilde \varrho}
        \end{pmatrix} u.
    \end{split}
\end{equation}
Now when $p\in (\ell^1_\nu)^2$ in \eqref{eq.zeroopbwd} corresponds to
$\left(\begin{smallmatrix}1\\0\end{smallmatrix}\right)$, we obtain the
  first row of $B(s)$; otherwise, when $p$ corresponds to
  $\left(\begin{smallmatrix}0\\1\end{smallmatrix}\right)$, we get the
    second row.

As in the forward flow case, computing the elements in
Lemma~\ref{lem:radii}, we end up proving the existence of the solution
$B(s) = \bar B(s) + \tilde B(s)$ with $\| \tilde B \| = \max _{i,j \in
  \{1,2\}} \| \tilde u _{ij}\| _{\ell ^1_\nu} \leq r _2$.

\subsection{Details of the Computer-Assisted Proofs for the Unperturbed Systems}

We present the results of three proofs. The first one for the periodic orbit and
the period is a nonlinear problem with a cubic term in the case of the
Van der Pol \eqref{vdp}.  The two other proofs, \eqref{eq.zeroopfwd}
and \eqref{eq.zeroopbwd}, are linear and they depend on the results
from the first proof.

We used the radii polynomial approach, see Appendix~\ref{app.radii},
for all the three computer-assisted proofs and a common $\nu=1.01$ for
the $\ell ^1_\nu$ space. When we encounter products, like in
$\vartheta$ and $\varrho$, which are convolutions in Chebyshev spaces,
we keep track of all the terms of the numerical and tail parts, see
Section \ref{sec:conv}. In particular, we have to manage quintic
convolution to get the norm of $D^2 K_0$ in the polynomials.

\begin{figure}[ht]
 \centering
 \includegraphics[scale=.55]{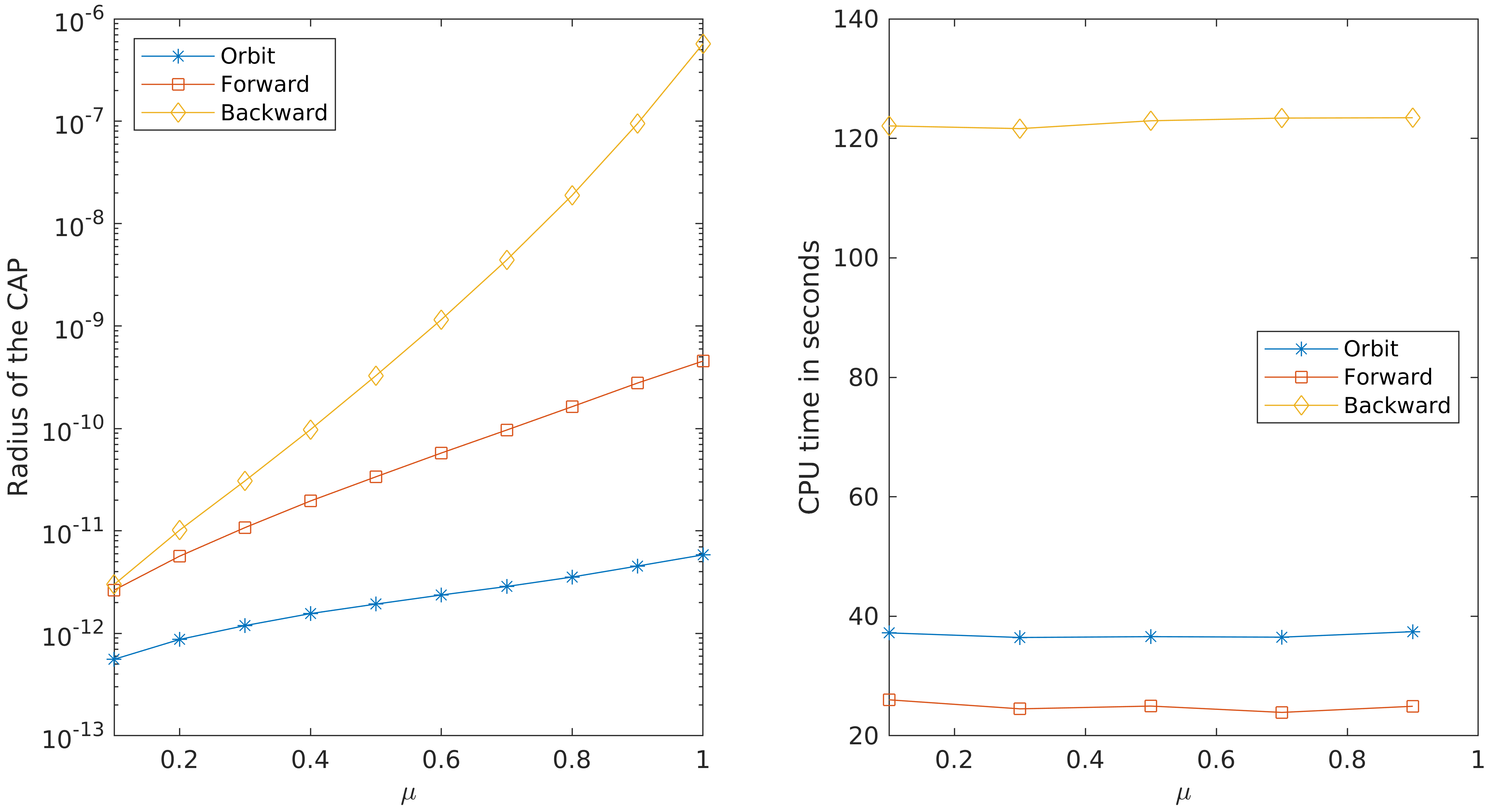}
 \caption{Radii of the computer-assisted proofs of the orbit, the
   forward, and the backward of \eqref{vdp} for different values of
   the parameter $\mu$.}
 \label{fig.ars}
\end{figure}

Figure~\ref{fig.ars} shows the final radii values and the
computational times of each proof, in particular, we observed that the
backward variational flow is computationally harder for larger $\mu$
because the periodic orbit becomes more attractive, see
Figure~\ref{fig.flows}.

\subsection{Bounds Computation and Optimization Steps}

Following Section~\ref{sec.compu-bounds}, we compute the bounds
required in Algorithm~\ref{alg.delay-cap}. We use different radii of
the computer-assisted proofs, that is, $r _0$ for the period and the
periodic orbit, $r _1$ for the forward flow, and $r _2$ for the
backward flow. Their values depend on the parameter $\mu$ in
\eqref{vdp}. For the triangle meshes we use the size of
$5000$. Overall the computation required around 3 days for each
parameter $\mu$ of the ODE \eqref{vdp}.

We have done two optimization processes. One in
Theorem~\ref{thm.vdp-sdde0} and another one in
Theorem~\ref{thm.vdp-sdde}. In both cases, the variables $a$ and
$\beta = (\beta _0, \beta _1,\beta _2)$ for the space $I _a \times
\mathcal{E} _ \beta$ in \eqref{space} have been restricted to the
domains
\begin{equation*}
 \begin{split}
  a, \beta _0 &\in (0, 0.1], \\
  \beta _1, \beta _2 &\in (0, 5], \\
  \varepsilon &\in (0, +\infty),
 \end{split}
\end{equation*}
and initial guesses $a = \beta _0 = 10^{-2}$, $\beta _1 = \beta _2 =
0.5$, and $\ep = 10^{-2}$.

\begin{figure}[ht]
 \centering
 \includegraphics[scale=.55]{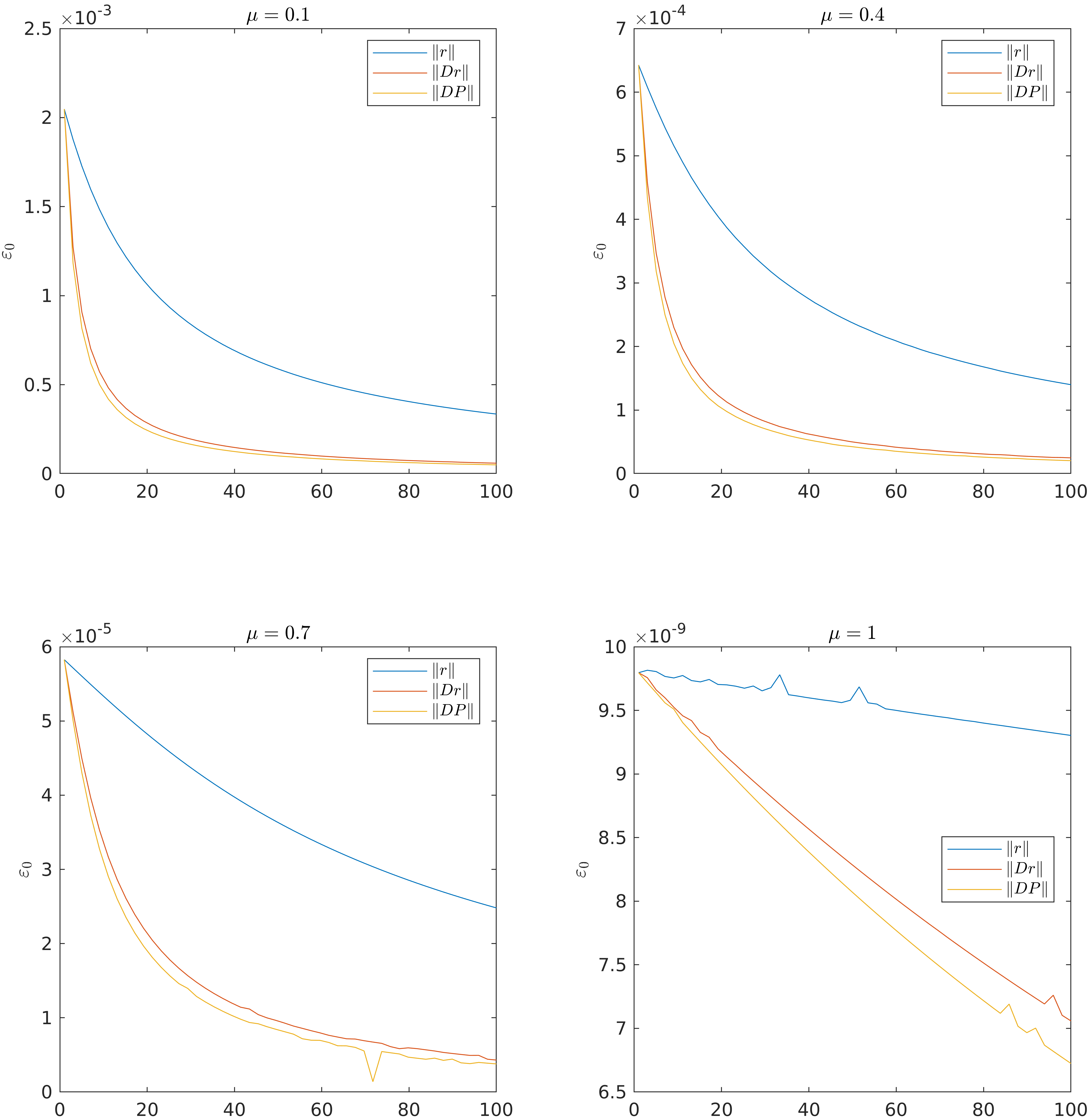}
 \caption{Value of the 
 $\varepsilon_0$ in \eqref{ineqs}
   for the Van der Pol equation \eqref{vdp}, moving either $\| r \|$,
   $\| Dr \|$, or $\| DP \|$ and setting the other two quatities equal
   to $1$.} \label{fig.optim_epsi}
\end{figure}

In Theorem~\ref{thm.vdp-sdde0} the objective function was just $\ep$
with $\|P\|$,
 $\|DP\|$, $\|r\|$, and $\|Dr\|$ equal to $1$. In
Figure~\ref{fig.optim_epsi} we used this procedure to illustrate, for
different values of the parameter $\mu$ in \eqref{vdp}, how sensitive the
numerical threshold $\ep _0$ is when either $\| r \|$, $\| D r \|$, or
$\| D P \|$ ranges in the $x$-axes of the plot and the other two
variables are set to $1$. Thus, from Figure~\ref{fig.optim_epsi}, we
observe that $\| DP \|$ and $\| Dr \|$ have similar, strong effects on
the $\ep_0$. On the other hand, $\|r\|$ presents less
influence in the $\ep_0$. In any case, that sensitivity will
always depend on the model itself and the inputs of the
Algorithm~\ref{alg.delay-cap}.

\begin{figure}[ht]
 \centering
 \includegraphics[scale=.55]{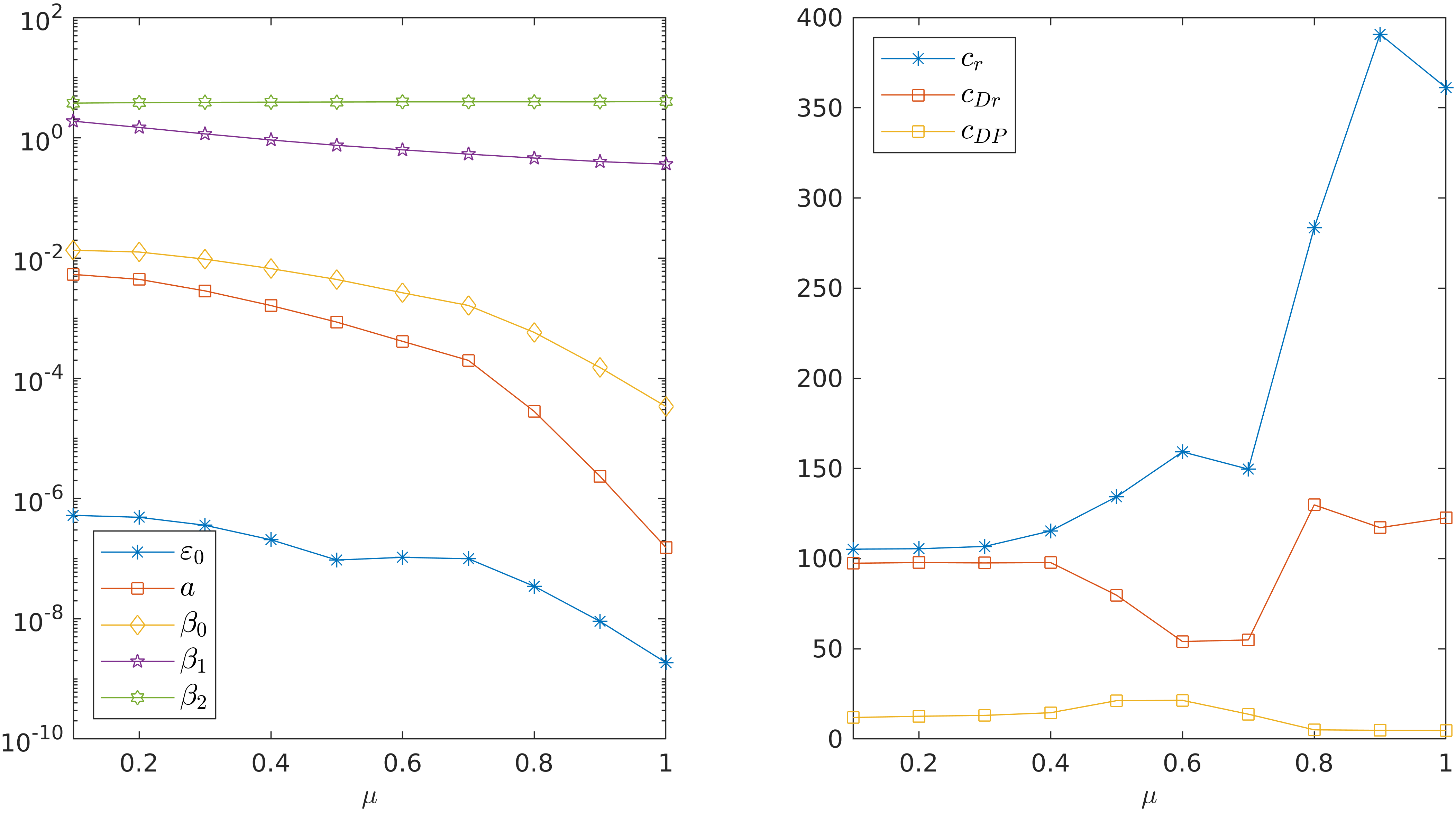}
 \caption{Values of the proof optimizing $\varepsilon$, $\| r \|$, $\|
   Dr \|$, and $\| DP \|$ for different values of the parameter $\mu$
   in \eqref{vdp} in Theorem~\ref{thm.vdp-sdde}.} \label{fig.capdelay}
\end{figure}

In Theorem~\ref{thm.vdp-sdde} the objective function is $\ep ^2 \| r
\| \| Dr \| \| DP \| ^2$ with $\| r \|$, $\| D r \|$, and $\| D P \|$
in the domain $[1,+\infty)$. Figure~\ref{fig.capdelay} shows the
  values $a$ and $\beta = (\beta _0,\beta _1,\beta _2 )$ for the space
  $I _a \times \mathcal{E} _ \beta$ in \eqref{space} and the numerical
  thresholds $\ep _0$, $c _{r}$, $c _{Dr}$, and $c _{DP}$ in
  Theorem~\ref{thm.vdp-sdde}.  Heuristically, the objective function
  was chosen considering that different operands have different
  scales. For instance, $\ep$ will be, in general, small and the rest
  will be large. Thus, we use multiplication instead of addition. We
  are putting $\ep ^2$ and $\|DP \|^2$ because we want them to be
  dominant during the optimization.  The reason is that the size of
  $\ep$ corresponds to the size of the perturbation, and it is not
  very interesting if the perturbation term $P$ is close to a constant
  vector.

We note that in our result, it is possible that none of $\ep$, $\|r
\|$, $\| Dr \|$, and $\| DP \|$ is maximized, since the optimization
processes are subject to various tolerances, the initial guess, etc.
Our goal is to provide some values of the parameters close to
optimized values so that the inequalities \eqref{ineqs} are satisfied.

Finally, after each of the numerical optimizations, we verify the
inequalities with the interval arithmetic polynomials we got before the
optimization to prove rigorously that they verify the inequalities
\eqref{ineqs}.  

\section{Conclusion} \label{sec:conclusion}

The theorems established here exploit well-known (but fairly state-of-the-art)
computer-assisted methods of proofs for ODEs,
however the use of these methods in the present work is fairly novel, as
the validated numerical computations are used to verify the hypotheses of  
a very singular perturbation theorem.  The statement of the theorem
hypothesizes the existence and other more quantitative properties  
of an isolated periodic orbit in a nonlinear system of ODEs. 
These hypotheses are notoriously difficult to verify in nonlinear ODEs, 
especially if the ODEs are far from any perturbative or asymptotic regime.  
The validated numerical methods allow us to pass from good numerical 
computations, to mathematically rigorous statements about the desired 
periodic solution.  Moreover, the high order spectral methods used in the 
present work provide enough control of the orbits that we can 
obtain, a-posteriori bounds for all the constants appearing in the 
hypotheses of the perturbation theorem.  Indeed, all of this can be made 
fairly automatic.

The result given here can be generalized easily to higher dimensional systems of 
ODEs, to the case when there are multiple forward or backward delays, either
state-dependent, distributed, or of other types. Of course the polynomials
considered here will need some modifications in other cases, 
but the constructions are not fundamentally different. 

Using the same framework, but
longer expressions for the polynomials, one should be able to study the case where
small delays are present as in electrodynamics, see \cite{Per}.
Indeed, applying the arguments developed here to relativistic perturbations 
of electrodynamics would be a fascinating future project.  
Our work used the \textsc{intlab} package, \cite{Ru99a} but other
packages providing rigorous bounds can be applied as well, such as,
\textsc{arb} in \cite{Johansson2017arb} and the CAPD library
\cite{MR4283203}.


\section*{Acknowledgment}
J.G. was supported by the Italian grant MIUR-PRIN 20178CJA2B ``New Frontiers of Celestial 
Mechanics: theory and Applications'', the Spanish grant PGC2018-100699-B-I00 (MCIU/AEI/FEDER, UE), 
and the Catalan grant 2017 SGR 1374. 
J.-P. L. was supported by an NSERC Discovery Grant. 
J.Y was supported by the National Science Foundation under Grant No. DMS-1929284 while in residence 
at the Institute for Computational and Experimental Research in Mathematics in Providence, 
RI, Fall 2021 - Spring 2022.
J.D. Mireles-James was partially supported by the National Science Foundation 
Grant DMS 1813501.  

\appendix
\section{Radii Polynomial Approach}
\label{app.radii}
There are many good references for
the radii polynomial approach, see for example, \cite{AJPJ,JPC,Mem},
we cite the one in \cite{Mem}, which suits our problem the best.
\begin{lemma}\label{lem:radii}
 Let $\bar x\in X$ and $r>0$, assume that $\mathcal G\colon X\to Y$ is
 Fr\'echet differentiable on the ball $B_r(\bar x)$. Consider bounded
 linear operators $A^\dagger\in B(X,Y)$ (approximation of $D \mathcal
 G(\bar x))$ and $A\in B(Y,X)$ (approximate inverse of $D \mathcal
 G(\bar x))$.

 Assume that $A$ is injective. Let $Y_0, Z_0, Z_1, Z_2\ge 0$ be bounds
 satisfying
\begin{align*}
\|A \mathcal G(\bar x)\|_X&\le Y_0\\
\|I-AA^\dagger\|_{B(X)}&\le Z_0\\
\|A(D \mathcal G\big(\bar x)-A^\dagger\big)\|_{B(X)}&\le Z_1\\
\|A\big(D \mathcal G(\bar x+z)-D \mathcal G(\bar x)\big)\|_{B(X)}&\le Z_2 r, \quad \forall z\in B_r(0) .
\end{align*}
Define the radii polynomial 
\begin{equation} \label{eq:radii_polynomial}
p(r)=Z_2 r^2+(Z_0+Z_1-1)r+Y_0 .
\end{equation}
If there exists $0<r_0\le r$, s.t. $p(r_0)<0$, then there exists a
unique $x\in B_{r_0}(\bar x)$ s.t. $\mathcal G(x)=0$.
\end{lemma}

Following \cite{JPC}, we briefly summarize the process to address
initial and boundary value problems. In both cases, we use a suitable
Banach space which is a Banach algebra under discrete convolutions.

\subsection{The Banach algebra space} \label{app:l1nu}
Let $\omega _k$ be the weights defined in terms of a parameter $\nu >
1$ and the Kronecker delta function $\delta _{i,j}$ as $\omega _k = (2
- \delta_{0,k})\nu ^k$. The ``ell one nu'' space is the normed-space
defined by
\begin{equation}
\label{eq.ell1nu}
    \ell ^1_\nu \bydef \biggl\{ \{b _k \} _{k \geq 0} \subset
    \mathbb{R} \colon \| b \| _ {\ell ^1_\nu} \bydef \sum _{k \geq
      0} \omega_k | b _k | < \infty \biggr \},
\end{equation}
In this space, we define a product (or convolution) as
\begin{equation}
 \label{eq.convo}
 (a \ast b) _k = \sum _{\substack{|k _1| + | k _2 | = k \\ k _1, k _2
     \in \mathbb{Z}}} a _{|k _1|} b _{|k _1|}
\end{equation}
for all $a = \{a _k\} _{k \geq 0}$ and $b = \{b _k\} _{k \geq 0}$ in
$\ell ^1_\nu$. The space is a Banach algebra under convolution
(e.g. see \cite{Lessard2018}),
\[
 \| a \ast b \| _ {\ell ^1_\nu}\leq \| a \| _ {\ell ^1_\nu} \| b \|_ {\ell ^1_\nu}.
\]

\subsubsection{The space of Chebyshev series}
\label{app.cheb}
The Chebyshev series of a function $x$ defined on $[-1,1]$ is an
expression of the form
\begin{equation}
\label{eq.cheb}
 x(t) = a _0 + 2 \sum _{k \geq 1} a _k T _k(t), \qquad a _0, a _k \in
 \R \text{ and } t \in [-1,1]
\end{equation}
where $T _k(t) = \cos (k \cos ^{-1}(t)) \in [-1,1]$ are the Chebyshev
polynomials, \cite{MR3012510}.

We consider the sequence of coefficients of \eqref{eq.cheb} in $\ell
^1_\nu$ for a parameter $\nu \geq 1$. Thus, the $\ell ^1_\nu$-norm of
the sequence is an upper bound of the $C^0$-norm of the function, i.e.
\[ 
 \| x \| _{C ^0} \leq \| x \| _{\ell ^1_1} \leq \| x \| _{\ell ^1_\nu} .
\]

The product of two Chebyshev series corresponds to the convolution
\eqref{eq.convo} of their coefficients. Note that the convolution
operator by $a$ in $\ell ^1_\nu$, i.e. $b \mapsto a \ast b$, can be
written as a combination of Toeplitz and Hankel (infinite) matrices.

From the numerical point of view, the Chebyshev series
\eqref{eq.convo} is truncated up to a finite number of coefficients
and its tail is estimated separately. Because of the aliasing
phenomenon, one needs to consider a zero-padding on the inputs to
obtain a correct output of convolution.

The Matlab code in Listing~\ref{lst.convo} codifies an arbitrary
number of convolutions using the Toeplitz and Hankel matrices. Note
that the code also works with the \textsc{intlab} package for interval
arithmetic. There are other approaches to compute rigorously the
convolution of Chebyshev series using the FFT,
e.g. \cite{Lessard2018}. This other approach is, in general, more
flexible when the memory resources is a constraint.

\begin{lstlisting}[caption = {Chebyshev convolutions using matrices.},
label={lst.convo}]
function a=convo(varargin) 
  m = nargin-1;
  nn = arrayfun(@(x) size(x{1},1),varargin);
  maxn = max(nn);
  a = [varargin{end}; zeros(maxn*(m+1)-nn(end),1)];
  for k = 1:m
   a = mat([varargin{k}; zeros(maxn*(m+1)-nn(k),1)])*a;
  end
end

function A=mat(a)
  m = length(a);
  A = toeplitz(a) + [zeros(m-1,1) hankel(a(2:end)); zeros(1,m)];
end
\end{lstlisting}

Another important part for truncated Chebyshev series is the
evaluation at a given $t \in [-1,1]$, there are various efficient
algorithms to compute it, such as the Clenshaw-Curtis
\cite{ClenshawC1960} or Laurent-Horner methods
\cite{AurentzH2020}. Using \textsc{intlab} we found that the direct
computation of the Chebyshev polynomials $T _k(t) = \cos (k \cos
^{-1}(t))$ and an explicit computation of the linear combination in
\eqref{eq.cheb} to be more accurate.

\subsection{Zero-finding problem for Initial Value Problems}
\label{app.ivp}
The first step for a rigorous proof of an initial value problem (IVP)
in $\R^d$, see \cite{JPC}, like
\begin{equation}
\label{eq.ivp}
 \dot x(t) = L f(x), \qquad x(-1) = p _0,
\end{equation}
consists in fixing the final time of integration that, by a temporal
change, encoded in the (known) parameter $L$, ensures $t \in [-1,1]$.
Note that \eqref{eq.ivp} is equivalent to
\begin{equation*}
 x(t) = x(-1) + L \int _{-1}^t f(x(s))\, ds, \qquad x(-1) = p _0.
\end{equation*}

If the orbit $x(t)$ is represented with Chebyshev coefficients $\{a _k
\} _{k \geq 0}$, as in \eqref{eq.cheb}, and we assume that the map $f$
admits a Chebyshev series representation, possibly by
polynomialization techniques \cite{MR4292532,MR3545977}, then $y(t) =
f(x(t))$ has computable coefficients $\{ f _k \} _{k \geq 0}$ in terms
of $\{a _k\}$. Thus, the equation $\mathcal{G}(a _k) =0$ for the
Lemma~\ref{lem:radii} corresponds to finding zero of the expression
\begin{equation*}
 \begin{cases}
   p _0  - a _0 - 2 \sum_{j \geq 1} (-1)^j a _j & k = 0, \\
   2 k a _{k} + L f _{k +1} - L f _{k-1} & k \geq 1.
 \end{cases}
\end{equation*}

\subsection{Zero-finding problem for the Boundary Value Problems}
\label{app.bvp}
Following \cite{JPC}, a periodic boundary value problem (BVP) of an
ODE consists in integrating an IVP, see Appendix~\ref{app.ivp}, taken
into consideration of $m$ boundary conditions given by the equation
\begin{equation*}
 \mathcal{H}(x(-1), x(1)) = 0
\end{equation*}
where $\mathcal{H} \colon \R ^{2n} \to \R ^m$ and $p _1=x(1)$ depend
on parameters $\theta \in \R ^m$, possibly including $L$ in
\eqref{eq.ivp} now as an unknown.

Assuming that $\mathcal{H}$ admits a Chebyshev representation series,
say $H(a _k, \theta)$, then the equation $\mathcal{G}(a _k, \theta)
=0$ for Lemma~\ref{lem:radii} corresponds to making zero the
expression
\begin{equation*}
 \begin{cases}
   H(a _k, \theta) & k =-1, \\
   p _1  - a _0 - 2 \sum_{j \geq 1} a _j & k = 0, \\
   2 k a _{k} + L f _{k +1} - L f _{k-1} & k \geq 1.
 \end{cases}
\end{equation*}

\section{The radii polynomial approach applied to the Van der Pol system}
The explicit bounds in Lemma~\ref{lem:radii} depends on the model
itself. Here we give one-by-one each of those bounds for the periodic
orbit of the Van der Pol case, \eqref{vdp}.  The forward and backward
variational flows are simpler since they are linear problems and, in
particular, $Z _2 = 0$.

First we need to consider the truncation up to Chebyshev coeefficient
order $n$, and bound the tails. After that, we can consider operator
$A$ and $A ^\dagger$ and we provide bounds for $Y _0$, $Z _0$, $Z _1$,
and $Z _2$.

\subsection{Projections and inclusions} \label{app.projinc}
Fixed $n \geq 0$. Define projection $\pi^n:\ell^1_\nu\to \R^{n+1}$ by
$(a _k) _{k \geq 0} \overset{\pi ^n}\mapsto (a_k)_{k=0}^n$. And,
similarly, define $\Pi^{(n)}\colon X_\nu \to \R^{2n+3}$,
\begin{equation*}
\Pi^{(n)}(a, b,L) \bydef (\pi^n a, \pi^n b, L), \qquad \forall\,
(a,b,L) \in X _\nu.
\end{equation*}

Let $\iota^n \colon \R^{n+1} \hookrightarrow \ell^1_\nu $ be the
inclusion defined by
\begin{equation*}
    (a _k) _{k = 0}^n \overset{\iota^n}\mapsto (a _0, a _1, \dotsc, a
_n, 0, 0, \dotsc).
\end{equation*}
Then similarly as in the projection case, we denote by
$\boldsymbol{\iota}^{(n)} \colon \R^{2n+3} \to X _\nu$ the mapping
\begin{equation*}
    \boldsymbol{\iota}^{(n)} (a , b , L) \bydef (\iota ^n a ,
    \iota^n b, L), \qquad \forall\, (a,b,L) \in (\R^{n+1} )^2\times
    \R.
\end{equation*}
For the periodic orbit, define a finite dimensional projection
$F^{(n)}\colon \R^{2n+3} \to \R^{2n+3}$ of $F$:
\begin{equation*}
    F^{(n)}(x) \bydef \Pi^{(n)}F(\boldsymbol{\iota}^{(n)} x).
\end{equation*}

\subsection{Choices of $A^\dagger$ and $A$ for the periodic orbit}
$A^\dagger$ is an approximation of the derivative of the zero-finding
function (at the numerical solution). We basically choose it as the
derivative of the truncated finite-dimensional zero-finding function
plus only the diagonal part for the tail.

More precisely, 
\begin{equation*}
    A^\dagger=
    \begin{pmatrix}
    A^\dagger_{1,a_1} & A^\dagger_{1,a_2} &
    A^\dagger_{1,L}\\
    A^\dagger_{2,a_1} & A^\dagger_{2,a_2} &
    A^\dagger_{2,L}\\
    A^\dagger_{3,a_1} & A^\dagger_{3,a_2} &
    A^\dagger_{3,L}\\
    \end{pmatrix} = 
    \begin{pmatrix}
    A^\dagger_{1,a_1} & 0 &
    A^\dagger_{1,L}\\
    0 & A^\dagger_{2,a_2} &
    A^\dagger_{2,L}\\
    A^\dagger_{3,a_1} & A^\dagger_{3,a_2} &
    0\\
    \end{pmatrix}.
\end{equation*}
For $b=(b_1, b_2, b_3)\in X _\nu$,
\begin{equation*}
    (A^\dagger b) _i = A^\dagger _{i, a_1} b _1 + A^\dagger _{i, a_2}
  b _2 + A^\dagger _{i, L} b _3 \in \ell _{\nu'}^1, \quad i = 1, 2, 3.
\end{equation*}
Note that for $i=1,2$, $(A^\dagger_{i,L}b_3)_k=(\partial_L
F_i^{(n)}(\bar x)b_3)_k$ for $k\le n$, $(A^\dagger_{i,L}b_3)_k=0$
otherwise. For $i=1, 2, 3$, $j=1, 2$:
\begin{equation*}
    (A^\dagger_{i,a_j} b_j)_k \bydef 
    \begin{cases}
    (D_{a_j^{(n)}}F_i^{(n)}b_j^{(n)})_k & k\le n,\\
    2k \delta_{ij} (b_j)_k & k>n,
    \end{cases}
\end{equation*}
where $\delta _{ij}$ is the Kronecker delta.

For $A$, we take numerical inverse of the derivative for the truncated
finite-dimensional zero-finding function and add the diagonal part for
the tail, which is smoothing.

Let $A^{(n)}$ be numerical inverse of $DF^{(n)}(\bar x)$,
\begin{equation*}
    A^{(n)} \bydef 
    \begin{pmatrix}
    A^{(n)}_{1,a_1} & A^{(n)}_{1,a_2} &
    A^{(n)}_{1,L}\\
    A^{(n)}_{2,a_1} & A^{(n)}_{2,a_2} &
    A^{(n)}_{2,L}\\
    A^{(n)}_{3,a_1} & A^{(n)}_{3,a_2} &
    A^{(n)}_{3,L}\\
    \end{pmatrix} \in \R^{(2n+3)\times(2n+3)}.
\end{equation*}

Define approximate inverse $A$ of $DF(\bar x)$ by
\begin{equation*}
    A \bydef 
    \begin{pmatrix}
    A_{1,a_1} & A_{1,a_2} &
    A_{1,L}\\
    A_{2,a_1} & A_{2,a_2} &
    A_{2,L}\\
    A_{3,a_1} & A_{3,a_2} &
    A_{3,L}\\
    \end{pmatrix} = 
    \begin{pmatrix}
    A_{1,a_1} & A_{1,a_2} &
    A^{(n)}_{1,L}\\
    A_{2,a_1} & A_{2,a_2} &
    A^{(n)}_{2,L}\\
    A_{3,a_1} & A_{3,a_2} &
    A^{(n)}_{3,L}\\
    \end{pmatrix}.
\end{equation*}

Where we abused the notation $A^{(n)}_{i,L}$ ($i=1,2$) for it to
denote an element in $\ell^1_\nu$ by adding a zero tail.  For $i=1, 2,
3$, $j=1, 2$, and $b = (b _1, b _2, b _3) \in X _{\nu '}$
\begin{equation*}
    (A_{i,a_j} b_j)_k=
    \begin{cases}
    (A_{i,a_j}^{(n)}b_j^{(n)})_k & k\le n, \\
    \frac{1}{2k} \delta _{ij} (b_j)_k & k>n.
    \end{cases}
\end{equation*}

\subsection{$Y$ bounds for periodic orbits}
With the zero-finding function and our choice of $A$, we are ready to
calculate the $Y$ bounds component-wise.
{\footnotesize
\begin{equation*}
    \begin{split}
        Y_1^{(0)} = &\left|\sum_{i=1}^2(A^{(n)}_{1,a_i}F^{(n)}_i(\bar x))_0+(A^{(n)}_{1,L}F^{(n)}_3(\bar x))_0\right|\\
    &+ 2\sum_{k \le n}
  \left|\sum_{i=1}^2(A^{(n)}_{1,a_i}F^{(n)}_i(\bar x))_k+(A^{(n)}_{1,L}F^{(n)}_3(\bar x))_k\right|\nu^k\\
   &+\frac{L}{n+1}|(a_2)_n|\nu^{n+1},
      \end{split}
\end{equation*}

\begin{equation*}
    \begin{split}
     Y_2^{(0)} = &\left|\sum_{i=1}^2(A^{(n)}_{2,a_i}F^{(n)}_i(\bar x))_0+(A^{(n)}_{2,L}F^{(n)}_3(\bar x))_0\right|\\
    &+ 2\sum_{k \le n} \left|\sum_{i=1}^2(A^{(n)}_{2,a_i}F^{(n)}_i(\bar x))_k+(A^{(n)}_{2,L}F^{(n)}_3(\bar x))_k\right|\nu^k\\
    &+\sum_{k=n+1}^{3n+1}\frac{L\mu}{k}\left|\sum_{k_1+k_2+k_3=k+1}(a_1)_{|k_1|}(a_1)_{|k_2|}(a_2)_{|k_3|}\right.\\
    &\phantom{AAAAAAAAAAAAAAAA}\left.-\sum_{k_1+k_2+k_3=k-1}(a_1)_{|k_1|}(a_1)_{|k_2|}(a_2)_{|k_3|}\right|\nu^k\\
     &\phantom{AA}+\frac{L}{n+1}|\mu (a_2)_n-(a_1)_n|\nu^{n+1},
    \end{split}
\end{equation*}

\begin{equation*}
    \begin{split}
        Y_3^{(0)} = &\left|\sum_{i=1}^2 A^{(n)}_{3,a_i}F^{(n)}_i(\bar x)+A^{(n)}_{3,L}F^{(n)}_3(\bar x)\right|.
    \end{split}
\end{equation*}}

\subsection{$Z$ bounds for periodic orbits}
\subsubsection{$Z^0$ bounds}

Let $B=I-AA^\dagger$, block-wise:
\begin{equation*}
    B \bydef 
    \begin{pmatrix}
     B _{1,a _1} & B _{1, a _2} & B _{1, L} \\
     B _{2,a _1} & B _{2, a _2} & B _{2, L} \\
     B _{3,a _1} & B _{3, a _2} & B _{3, L}
    \end{pmatrix}.
\end{equation*}
Define the norm
\[
\begin{split}
\| B _{j,a _i} \| _{B(\ell _\nu^1, \ell _\nu^1)}^* \bydef &\sup_{m\in\N}\frac{1}{\omega_m} \sum_{k\in \N}|(B _{j, a _i}) _{k, m}|\omega_k\\ =&\max_{0\le m \le n}\frac{1}{\omega_m} \sum_{0\le k \le n}|(B _{j, a _i}) _{k, m}|\omega_k, \qquad \forall\, i, j\in \{1,2\}.
\end{split}
\]
Then, we can take
\begin{equation*}
\begin{split}
    Z _j ^{(0)} &\bydef  \| B _{j,a _1} \| _{B(\ell _\nu^1, \ell _\nu^1)}^* + \| B _{j,a _2} \| _{B(\ell _\nu^1, \ell _\nu^1)}^* +  \| B _{j,L} \| _\nu, \qquad \forall\, j \in \{1,2\}, \\ 
    Z _3 ^{(0)} &\bydef \| B _{3,a_1} \|_\nu^\infty + \| B _{3,a_2} \|_\nu^\infty + |B _{3,L}| ,
\end{split}
\end{equation*}

\subsubsection{$Z^1$ bounds}
let $z_j=([DF(\bar x)-A^\dagger]d)_j$. Let $P^{n}=\iota^n\circ\pi^n \colon \ell^1_\nu\to\ell^1_\nu$, 
\[
 P^n\colon (a _k) _{k \ge 0} \mapsto (a _0, a _1, \dotsc, a _n, 0, 0, \dotsc),
 \]
 and $P^I=Id-P^{n} \colon \ell^1_\nu\to\ell^1_\nu$,
 \[
 P^I\colon (a _k) _{k \ge 0} \mapsto (0, 0, \dotsc,0, a _{n+1},a _{n+2}, \dotsc).
 \]
Define operator $T^I\bydef T-P^n T P^n$.

\begin{equation*}
 z_1=M P^I d_1+\bar L\cdot T^I d_2+d_3\cdot P^I T \bar a_2.
\end{equation*}

Notice that
\begin{equation*}
    Dc_2(a)(d_1, d_2)=\mu d_2-2\mu(d_1*a_1*a_2)-\mu(a_1*a_1*d_2)-d_1,
\end{equation*}
\begin{equation*}
\begin{split}
    z_2=&M P^I d_2+\bar L\cdot T^I (\mu d_2-d_1)+d_3\cdot P^I T c_2( \bar a)\\
&-2\mu\bar LP^n T(P^Id_1*\bar a_1*\bar a_2)-2\mu\bar LP^I T(d_1*\bar a_1*\bar a_2), \\
&-\mu\bar LP^n T(\bar a_1*\bar a_1*P^Id_2)-\mu\bar LP^I T(\bar a_1*\bar a_1*d_2),
\end{split}
\end{equation*}

\begin{equation*}
     z _3 = 2 \sum _{k> n} (-1)^k (d_1)_k.
\end{equation*}

Note that
\begin{equation*}
    \|MP^I\|_{B(\ell^1_\nu, \ell^1_\nu)}\le \frac{1}{\omega_{n+1}},
\end{equation*}
\begin{equation*}
    \|T^I\|_{B(\ell^1_\nu, \ell^1_\nu)}\le \frac{1}{\nu}+\nu,
\end{equation*}
\begin{equation*}
    \begin{split}
        \|A_{1,a_1}z_1\|_\nu\le &\frac{1}{\omega_{n+1}}\left(\|(A_{1,a_1}^{(n)})_{.,0}\|_\nu+\bar L\|(A_{1,a_1}^{(n)})_{.,n}\|_\nu\right)r\\
        &+\frac{1}{2(n+1)}\left(\bar L (\frac{1}{\nu}+\nu)+|(\bar a_2)_n|\omega_{n+1}\right)r,
    \end{split}
\end{equation*}

\begin{equation*}
     \|A_{2,a_1}z_1\|_\nu\le \frac{1}{\omega_{n+1}}\left(\|(A_{2,a_1}^{(n)})_{.,0}\|_\nu+\bar L\|(A_{2,a_1}^{(n)})_{.,n}\|_\nu\right)r,
\end{equation*}

\begin{equation*}
     |A_{3,a_1}z_1|\le \frac{1}{\omega_{n+1}}\left(|(A_{3,a_1}^{(n)})_0|+\bar L|(A_{3,a_1}^{(n)})_n|\right)r.
\end{equation*}

To bound $z_2$, we need the following lemma.

\begin{lemma}
Let $N\in \N$, and $\bar \alpha=(\bar \alpha_0, \bar \alpha_1, \bar
\alpha_2,\dotsc,\bar \alpha_N,0, 0, \dotsc)\in \ell^1_\nu$. For $0\le
k\le n+1$, define $\hat{\ell}^k_{\bar \alpha}\in\ell^\infty_\nu$ by
\begin{equation*}
    \hat{\ell}^k_{\bar \alpha}(h)\bydef (\bar \alpha*P^I h)_k=\sum_{\substack{k_1+k_2=k \\k_1, k_2\in \Z}}\bar \alpha_{|k_1|}(P^I h)_{|k_2|}.
\end{equation*}
Then, 
\begin{equation*}
    \|\hat{\ell}^k_{\bar \alpha}\|_\nu^\infty\le \Psi_k(\bar \alpha)\bydef \max_{n<j\le k+N}\left(\frac{\left|\bar\alpha_{|k-j|}+\bar\alpha_{k+j}\right|}{2\nu^j}\right).
\end{equation*}
\end{lemma}
In our case, we will let $N=2n$.
{\footnotesize
\begin{equation*}
    \begin{split}
        \|A_{1,a_2}z_2\|_\nu\le &\frac{1}{\omega_{n+1}}\left(\|(A_{1,a_2}^{(n)})_{.,0}\|_\nu+(\mu+1)\bar L\|(A_{1,a_2}^{(n)})_{.,n}\|_\nu\right)r\\
        &+2\mu \bar L \sum^n_{i=0}\sum^n_{j=1}|(A_{1,a_2}^{(n)})_{i,j}|\left(\Psi_{j-1}(\bar a_1*\bar a_2)+\Psi_{j+1}(\bar a_1*\bar a_2)\right)\omega_i r\\
        &+\mu \bar L \sum^n_{i=0}\sum^n_{j=1}|(A_{1,a_2}^{(n)})_{i,j}|\left(\Psi_{j-1}(\bar a_1*\bar a_1)+\Psi_{j+1}(\bar a_1*\bar a_1)\right)\omega_i r,
        \end{split}
\end{equation*}

\begin{equation*}
    \begin{split}
        \|A_{2,a_2}z_2\|_\nu\le &\frac{1}{\omega_{n+1}}\left(\|(A_{2,a_2}^{(n)})_{.,0}\|_\nu+(\mu+1)\bar L\|(A_{2,a_2}^{(n)})_{.,n}\|_\nu\right)r\\
        &+2\mu \bar L \sum^n_{i=0}\sum^n_{j=1}|(A_{2,a_2}^{(n)})_{i,j}|\left(\Psi_{j-1}(\bar a_1*\bar a_2)+\Psi_{j+1}(\bar a_1*\bar a_2)\right)\omega_i r\\
        &+\mu \bar L \sum^n_{i=0}\sum^n_{j=1}|(A_{2,a_2}^{(n)})_{i,j}|\left(\Psi_{j-1}(\bar a_1*\bar a_1)+\Psi_{j+1}(\bar a_1*\bar a_1)\right)\omega_i r\\
        &+\frac{1}{2(n+1)}\left(|(\mu\bar a_2-\bar a_1)_n|\omega_{n+1}+\bar L (\mu+1)(\frac{1}{\nu}+\nu)\right)r\\
        &+\frac{\mu\|\bar a_1\|_\nu}{2(n+1)}(\frac{1}{\nu}+\nu)\left(\|\bar a_1\|_\nu\|\bar a_2\|_\nu+\bar L (\|\bar a_2\|_\nu+\|\bar a_1\|_\nu)\right)r,
    \end{split}
\end{equation*}

\begin{equation*}
    \begin{split}
        \|A_{3,a_2}z_2\|_\nu\le &\frac{1}{\omega_{n+1}}\left(|(A_{3,a_2}^{(n)})_0|+(\mu+1)\bar L|(A_{3,a_2}^{(n)})_n|\right)r\\
        &+2\mu \bar L \sum^n_{i=1}|(A_{3,a_2}^{(n)})_i|\left(\Psi_{i-1}(\bar a_1*\bar a_2)+\Psi_{i+1}(\bar a_1*\bar a_2)\right) r\\
        &+\mu \bar L \sum^n_{i=1}|(A_{3,a_2}^{(n)})_i|\left(\Psi_{i-1}(\bar a_1*\bar a_1)+\Psi_{i+1}(\bar a_1*\bar a_1)\right) r.
        \end{split}
\end{equation*}}
Since $|z_3|\le \frac{2r}{\omega_{n+1}}$,
\begin{equation*}
    \|A_{1,L}z_3\|_\nu\le\frac{2}{\omega_{n+1}}\|A_{1,L}^{(n)}\|_\nu r,
\end{equation*}
\begin{equation*}
    \|A_{2,L}z_3\|_\nu\le\frac{2}{\omega_{n+1}}\|A_{2,L}^{(n)}\|_\nu r,
\end{equation*}
\begin{equation*}
    |A_{3,L}z_3|\le\frac{2}{\omega_{n+1}}|A_{3,L}^{(n)}|r.
\end{equation*}
Now we let 
{\footnotesize
\begin{equation*}
    \begin{split}
        Z^{(1)}_1\bydef &\frac{\left(\sum^2_{k=1}\|(A_{1,a_k}^{(n)})_{.,0}\|_\nu+\bar L\left(\|(A_{1,a_1}^{(n)})_{.,n}\|_\nu+(\mu+1)\|(A_{1,a_2}^{(n)})_{.,n}\|_\nu\right) + 2\|A_{1,L}^{(n)}\|_\nu\right)}{\omega_{n+1}}\\
        &+\frac{1}{2(n+1)}\left(\bar L (\frac{1}{\nu}+\nu)+|(\bar a_2)_n|\omega_{n+1}\right)\\
        &+2\mu \bar L \sum^n_{i=0}\sum^n_{j=1}|(A_{1,a_2}^{(n)})_{i,j}|\left(\Psi_{j-1}(\bar a_1*\bar a_2)+\Psi_{j+1}(\bar a_1*\bar a_2)\right)\omega_i\\
        &+\mu \bar L \sum^n_{i=0}\sum^n_{j=1}|(A_{1,a_2}^{(n)})_{i,j}|\left(\Psi_{j-1}(\bar a_1*\bar a_1)+\Psi_{j+1}(\bar a_1*\bar a_1)\right)\omega_i ,
    \end{split}
\end{equation*}

\begin{equation*}
    \begin{split}
        Z^{(1)}_2\bydef & \frac{\sum^2_{k=1}\|(A_{2,a_k}^{(n)})_{.,0}\|_\nu+\bar L\left(\|(A_{2,a_1}^{(n)})_{.,n})\|_\nu+(\mu+1)\|(A_{2,a_2}^{(n)})_{.,n}\|_\nu\right) + 2\|A_{2,L}^{(n)}\|_\nu}{\omega_{n+1}}\\
         &+2\mu \bar L \sum^n_{i=0}\sum^n_{j=1}|(A_{2,a_2}^{(n)})_{i,j}|\left(\Psi_{j-1}(\bar a_1*\bar a_2)+\Psi_{j+1}(\bar a_1*\bar a_2)\right)\omega_i \\
        &+\mu \bar L \sum^n_{i=0}\sum^n_{j=1}|(A_{2,a_2}^{(n)})_{i,j}|\left(\Psi_{j-1}(\bar a_1*\bar a_1)+\Psi_{j+1}(\bar a_1*\bar a_1)\right)\omega_i \\
        &+\frac{1}{2(n+1)}\left(|(\mu\bar a_2-\bar a_1)_n|\omega_{n+1}+\bar L (\mu+1)(\frac{1}{\nu}+\nu)\right)\\
        &+\frac{\mu}{2(n+1)}(\frac{1}{\nu}+\nu)\left(\|\bar a _1 * \bar a_1 * \bar a_2\|_\nu+\bar L (2\|\bar a _1 *\bar a_2\|_\nu+\|\bar a _1 * \bar a_1\|_\nu)\right),
    \end{split}
\end{equation*}

\begin{equation*}
    \begin{split}
        Z^{(1)}_3\bydef&\frac{1}{\omega_{n+1}}\left(\sum^2_{k=1}|(A_{3,a_k}^{(n)})_0|+\bar L\left(|(A_{3,a_1}^{(n)})_n|+(\mu+1)|(A_{3,a_2}^{(n)})_n|\right) + 2 |A_{3,L}^{(n)}|\right)\\
      &+2\mu \bar L \sum^n_{i=1}|(A_{3,a_2}^{(n)})_i|\left(\Psi_{i-1}(\bar a_1*\bar a_2)+\Psi_{i+1}(\bar a_1*\bar a_2)\right)\\
        &+\mu \bar L \sum^n_{i=1}|(A_{3,a_2}^{(n)})_i|\left(\Psi_{i-1}(\bar a_1*\bar a_1)+\Psi_{i+1}(\bar a_1*\bar a_1)\right).
        \end{split}
\end{equation*}}

\subsubsection{Z\textsuperscript{(2)} bound}
Now let $\zeta_j=([DF(\bar x +b)-DF(\bar x)]d)_j$. Then, $(\zeta_i)_0=0$ for $i=1,2$, and $\zeta_3=0$.

For $k\ge 1$, 
\begin{equation*}
\begin{split}
    \zeta_1&=d_3\cdot T b_2+b_3\cdot T d_2, \\
    \zeta_2&=[(\bar L+b_3)\cdot T Dc_2(\bar a +(b_1,b_2))-\bar L\cdot T D c_2(\bar a)](d_1, d_2)\\
    & \qquad +d_3\cdot T [c_2(\bar a +(b_1,b_2))-c_2(\bar a))].
\end{split}
\end{equation*}
Since
\begin{equation*}
    \|T\|_{B(\ell^1_\nu, \ell^1_\nu)}\le \frac{1}{\nu}+\nu,
\end{equation*}
then
\begin{equation*}
    \|\zeta_1\|_\nu\le 2\big(\frac{1}{\nu}+\nu\big)r^2
\end{equation*}
\begin{equation*}
\begin{split}
    \|\zeta_2\|_\nu\le& \big(\frac{1}{\nu}+\nu\big)\bigg((\mu+1)+2\mu(\|\bar a_1\|_\nu+r)(\|\bar a_2\|_\nu+r)+\mu(\|\bar a_1\|_\nu+r)^2\\
    &\phantom{AAAAAAA}+\mu\bar L\big( 2(\|\bar a_1\|_\nu+\|\bar a_2\|_\nu+r)+2\|\bar a_1\|_\nu+r\big) \bigg)r^2\\
    &+\big(\frac{1}{\nu}+\nu\big)\bigg((\mu+1)+\mu(\|\bar a_1\|_\nu^2+2\|\bar a_1\|_\nu\|\bar a_2\|_\nu+\|\bar a_2\|_\nu r+2\|\bar a_1\|_\nu r+ r^2)\bigg)r^2\\
    \le & \big(\frac{1}{\nu}+\nu\big)\bigg(\mu\big(4\|\bar a_1\|_\nu\|\bar a_2\|_\nu+2\|\bar a_1\|_\nu^2+6\|\bar a_1\|_\nu r+3\|\bar a_2\|_\nu r+4r^2\big)\\
    &\phantom{AAAAAAA}+\mu\bar L\big( 4\|\bar a_1\|_\nu+2\|\bar a_2\|_\nu+3r\big)+2(\mu+1) \bigg)r^2.
\end{split}
\end{equation*}

\begin{lemma}
For $i,j=1,2$,
\begin{equation*}
    \| A _{i,a _j} \| _{B(\ell _\nu^1, \ell _\nu^1)}^*\le \max \left\{\| A _{i,a _j}^{(n)} \| ^*,\frac{1}{2(n+1)}\delta_{i,j}\right\},
\end{equation*}
where $\| A _{i,a _j}^{(n)} \| ^*\bydef\max\limits_{0\le m \le n}\frac{1}{\omega_m} \sum\limits_{0\le k \le n}|(A _{i,a _j}^{(n)}) _{k, m}|\omega_k$.
\end{lemma}

Therefore we have
\begin{equation*}
    \|(A\zeta)_i\|_\nu\le  \sum^2_{j=1}\| A _{i,a _j} \| _{B(\ell _\nu^1, \ell _\nu^1)}^*\|\zeta_j\|_\nu, \qquad
     |(A\zeta)_3|\le  \sum^2_{j=1}\| A _{i,a _j} \| _\nu^\infty\|\zeta_j\|_\nu.
\end{equation*}
Define 
\begin{multline*}
    \Xi \bydef \mu\big(4\|\bar a_1*\bar a_2\|_\nu+2\|\bar a_1 * \bar a_1\|_\nu + 6\|\bar a_1\|_\nu r+3\|\bar a_2\|_\nu r+4r^2\big) \\
    +\mu\bar L\big( 4\|\bar a_1\|_\nu+2\|\bar a_2\|_\nu+3r\big)+2(\mu+1).
\end{multline*}
Then 
\begin{equation*}
\begin{split}
    Z^{(2)}_1&=\big(\frac{1}{\nu}+\nu\big)\| A _{1,a _2}^{(n)} \| ^*\Xi +2\big(\frac{1}{\nu}+\nu\big)\max \left\{\| A _{1,a _1}^{(n)} \| ^*,\frac{1}{2(n+1)}\right\}, \\
    Z^{(2)}_2&=\big(\frac{1}{\nu}+\nu\big)\max \left\{\| A _{2,a _2}^{(n)} \| ^*,\frac{1}{2(n+1)}\right\}\Xi +2\big(\frac{1}{\nu}+\nu\big)\| A _{2,a _1}^{(n)} \| ^*, \\ 
    Z^{(2)}_3&=\big(\frac{1}{\nu}+\nu\big)\| A _{3,a _2}^{(n)} \| _\nu^\infty\Xi +2\big(\frac{1}{\nu}+\nu\big)\| A _{3,a _1} ^{(n)}\| _\nu^\infty.
\end{split}
\end{equation*}

\section{Expansion of the polynomials in monomials} \label{app:monopoly}

The polynomials in Section~\ref{sec.polys} are used in the
optimization step of the Algorithm~\ref{alg.delay-cap}. However,
during the optimization is better to expand them in monomials to
provide sharper bounds.

The polynomials \eqref{q}, \eqref{p0}, and \eqref{p1} are already
expanded in monomials of the optimization variables. On the other
hand, The other polynomials \eqref{p2}, \eqref{mu1}, and \eqref{mu2}
have different monomial expansions depending on the optimization
variables.

\subsection{Monomial expansion for fixed perturbations}
In this case, the variables to be optimize are $\ep$, $ a$, $\beta
_0$, $\beta _1$, and $\beta _ 2$. Then
\begin{multline*}
 P _2(\ep, a, \beta _0, \beta _1, \beta _2) = 
 \frac{\| DP \| \| Dr \|}{\omega _0 } \ep a \beta _1^2 + 
 \frac{2 \| DK _0 \| \| DP \| \| Dr \|}{\omega _0 } \ep a \beta _1 + \\ 
 \frac{\| DK _0 \|^2 \| DP \| \| Dr \|}{\omega _0 } \ep a + 
 \| DP \| \| Dr \| \ep \beta _1 ^2 + \\
 \| DP \| \left(2 \| DK _0 \| \| Dr \|+\frac{1}{\omega _0 }\right) \ep \beta _1 +
 \frac{\| DK _0 \| \| DP \| (\| DK _0 \| \| Dr \|
   \omega _0 +1)}{\omega _0 } \ep + \\
 \frac{1}{\omega _0 } a \beta _2 +
 \frac{\| D ^2K _0 \|}{\omega _0 } a + 
 \frac{\| D ^3f \|}{\omega _0 } \beta _0^2 \beta _1 +
 \frac{\| D ^3f \| \| DK _0 \|}{\omega _0 } \beta _0^2 + \\
 \frac{2 \| D ^2f \|}{\omega _0 } \beta _0 \beta _1 +
 \frac{\| D ^2f \circ K _0 \| \| DK _0 \|}{\omega _0 } \beta _0 +
 \frac{\| D f \circ K _0 \|}{\omega _0 } \beta _1 .
\end{multline*}

\begin{multline*}
 \mu _1(\ep, a, \beta _0, \beta _1, \beta _2) = \frac{\| \Pi _{\theta _0}^\top \|C _{1,1}}{| DK _0(\theta _0) |} \biggl(
 \| DP \| \| Dr \|  \ep a \beta _1 + \\
 \| DK _0 \| \| DP \| \| Dr \| \ep a + 
  \| DP \|  (\| Dr \|
   \omega _0 +\| r \|) \ep \beta _1 + \\
  \| DP \|  (\| DK _0 \| (\| Dr \|
   \omega _0 +\| r \|)+1) \ep \biggr) + \\
 \frac{\| \Pi _{\theta _0}^\top \| C _{2,1} }{| DK _0(\theta _0) |} a +
 \frac{\| \Pi _{\theta _0}^\top \|
   (C _{1,1} \| D ^2f \|+C _{2,1})}{| DK _0(\theta _0) |}\beta _0 .
\end{multline*}
Let us define
\begin{equation} \label{eq.gammadelta} 
\begin{split}
 \gamma &\bydef \frac{C _{1,2}C _{1,1}
   \| DK _0 \| \| \Pi _{\theta _0}^\top \|}{| DK _0(\theta _0) | \omega _0 }+\frac{C _{1,2}}{ \omega _0 }+C _{1,1} C _{1,3}  M \\
  \delta &\bydef \frac{C _{1,2} C _{2,1} \| DK _0 \| \| \Pi _{\theta _0}^\top \|}{| DK _0(\theta _0) | \omega _0 }+C _{1,3} C _{2,1} M +\frac{C _{2,2}}{\omega _0},
\end{split}
\end{equation}
then 
\begin{multline*}
 \mu _2(\ep, a, \beta _0, \beta _1, \beta _2) = 
 \| DP \| \| Dr \| \gamma \ep a \beta _1 + \\
 \| DK _0 \| \| DP \| \| Dr \|  \gamma \ep a + 
 \| DP \| (\| Dr \| \omega _0 +\| r \|)  \gamma \ep \beta _1  + \\
 \| DP \| (\| DK _0 \| (\| Dr \| \omega _0 +\| r \|)+1) \gamma \ep +
 \delta a + (\| D ^2f \|\gamma + \delta)\beta _0 .
\end{multline*}

\subsection{Monomial expansion for class of perturbations}
In this case, the variables are $\ep$, $ a$, $\beta _0$, $\beta _1$,
$\beta _ 2$, $\| DP \|$, $\| Dr \|$, and $\|r \|$. Using
\eqref{eq.gammadelta} for \eqref{mu2},
\begin{multline*}
 P _2(\ep, a, \beta _0, \beta _1, \beta _2, \| DP \|, \|Dr \| ) = 
 \frac{1}{\omega _0 } \ep a \beta _1^2 \|DP \| \| Dr \| + \\
 \frac{2 \| DK _0 \|}{\omega _0 } \ep a \beta _1  \|DP \| \| Dr \| +
 \frac{\| DK _0 \|^2}{\omega _0 } \ep a  \|DP \| \| Dr \| + \\
 \ep \beta _1 ^2  \|DP \| \| Dr \| +
 2 \| DK _0 \| \ep \beta _1  \|DP \| \| Dr \| +
 \frac{1}{\omega _0 } \ep \beta _1  \|DP \|  +\\
 \| DK _0 \|^2 \ep  \|DP \| \| Dr \| +
 \frac{\| DK _0 \|}{\omega _0 } \ep  \| DP \|  +
 \frac{1}{\omega _0 } a \beta _2 + \\
 \frac{\| D ^2K \|}{\omega _0 } a +
 \frac{\| D ^3f \|}{\omega _0 } \beta _0^2 \beta _1 +
 \frac{\| D ^3f \|
   \| DK _0 \|}{\omega _0 } \beta _0^2 + \\
 \frac{2
   \| D ^2f \|}{\omega _0 } \beta _0 \beta _ +
 \frac{\| D ^2f \circ K _0 \|
   \| DK _0 \|}{\omega _0 } \beta _0 +
 \frac{\| D f \circ K _0 \|}{\omega _0 } \beta _1 .
\end{multline*}

\begin{multline*}
 \mu _1(\ep, a, \beta _0, \beta _1, \beta _2, \| DP \|, \|Dr \|, \| r \| ) = 
 \frac{\| \Pi _{\theta _0}^\top \| C _{1,1}}{| DK _0(\theta _0) |} \biggl(\ep a \beta _1  \|DP \| \| Dr \| + \\
 \| DK _0 \| \ep a \|DP \| \| Dr \| + 
 \omega _0 \ep \beta _1 \|DP \| \| Dr \| + 
 \ep \beta _1 \|DP \| \| r \| + \\
 \| DK _0 \|\omega _0  \ep  \|DP \| \| Dr \| + 
 \| DK _0 \| \ep  \|DP \| \| r \| + 
 \ep  \|DP \| \biggr) + \\
 \frac{
   \| \Pi _{\theta _0}^\top \|C _{2,1}}{| DK _0(\theta _0) |} a + 
 \frac{\| \Pi _{\theta _0}^\top \|
   (C _{1,1} \| D ^2f \|+C _{2,1})}{| DK _0(\theta _0) |} \beta _0,
\end{multline*}
and
\begin{multline*}
 \mu _2(\ep, a, \beta _0, \beta _1, \beta _2, \| DP \|, \|Dr \|, \| r \|) = 
 \gamma \ep a \beta _2  \|DP \| \| Dr \| + \\
 \| DK _0 \| \gamma \ep a \|DP \| \| Dr \| + 
 \omega _0 \gamma\ep \beta _2 \|DP \| \| Dr \| + \\
 \gamma \ep \beta _2 \|DP \| \| r \| +
 \| DK _0 \| \omega _0 \gamma \ep  \|DP \| \| Dr \| + \\
 \| DK _0 \| \gamma \ep  \|DP \| \| r \| + 
 \gamma \ep  \|DP \| + 
 \delta a + 
 (\| D ^2f \| \gamma + \delta) \beta _0 .
\end{multline*}

\bibliographystyle{alpha}
\bibliography{ref}

\newcommand{\etalchar}[1]{$^{#1}$}
\begin{thebibliography}{DvGVLW95}

\bibitem[AH20]{AurentzH2020}
Jared~L. Aurentz and Behnam Hashemi.
\newblock The {L}aurent-{H}orner method for validated evaluation of {C}hebyshev
  expansions.
\newblock {\em Appl. Math. Lett.}, 102:106113, 5, 2020.

\bibitem[CC60]{ClenshawC1960}
C.~W. Clenshaw and A.~R. Curtis.
\newblock A method for numerical integration on an automatic computer.
\newblock {\em Numer. Math.}, 2:197--205, 1960.

\bibitem[CCdlL20]{MR4112213}
Alfonso Casal, Livia Corsi, and Rafael de~la Llave.
\newblock Expansions in the delay of quasi-periodic solutions for state
  dependent delay equations.
\newblock {\em J. Phys. A}, 53(23):235202, 20, 2020.

\bibitem[CdlL20]{MR4066033}
Hongyu Cheng and Rafael de~la Llave.
\newblock Stable manifolds to bounded solutions in possibly ill-posed {PDE}s.
\newblock {\em J. Differential Equations}, 268(8):4830--4899, 2020.

\bibitem[CFdlL03a]{MR1976079}
X.~Cabr{\'e}, E.~Fontich, and R.~de~la Llave.
\newblock The parameterization method for invariant manifolds. {I}. {M}anifolds
  associated to non-resonant subspaces.
\newblock {\em Indiana Univ. Math. J.}, 52(2):283--328, 2003.

\bibitem[CFdlL03b]{MR1976080}
X.~Cabr{\'e}, E.~Fontich, and R.~de~la Llave.
\newblock The parameterization method for invariant manifolds. {II}.
  {R}egularity with respect to parameters.
\newblock {\em Indiana Univ. Math. J.}, 52(2):329--360, 2003.

\bibitem[CFdlL05]{MR2177465}
X.~Cabr{\'e}, E.~Fontich, and R.~de~la Llave.
\newblock The parameterization method for invariant manifolds. {III}.
  {O}verview and applications.
\newblock {\em J. Differential Equations}, 218(2):444--515, 2005.

\bibitem[CGL18]{MR3749257}
Roberto Castelli, Marcio Gameiro, and Jean-Philippe Lessard.
\newblock Rigorous numerics for ill-posed {PDE}s: periodic orbits in the
  {B}oussinesq equation.
\newblock {\em Arch. Ration. Mech. Anal.}, 228(1):129--157, 2018.

\bibitem[Chu21]{kevinIntegrator}
Kevin Church.
\newblock Validated integration of differential equations with state-dependent
  delay.
\newblock {\em (submitted)}, 2021.

\bibitem[CLMJ15]{MR3436565}
Roberto Castelli, Jean-Philippe Lessard, and Jason~D. Mireles~James.
\newblock Analytic enclosure of the fundamental matrix solution.
\newblock {\em Appl. Math.}, 60(6):617--636, 2015.

\bibitem[DL16]{MR3456815}
Jayme De~Luca.
\newblock Equations of motion for variational electrodynamics.
\newblock {\em J. Differential Equations}, 260(7):5816--5833, 2016.

\bibitem[DLGHP10]{MR2639911}
Jayme De~Luca, Nicola Guglielmi, Tony Humphries, and Antonio Politi.
\newblock Electromagnetic two-body problem: recurrent dynamics in the presence
  of state-dependent delay.
\newblock {\em J. Phys. A}, 43(20):205103, 20, 2010.

\bibitem[DLHR12]{MR2912694}
Jayme De~Luca, A.~R. Humphries, and Savio~B. Rodrigues.
\newblock Finite element boundary value integration of {W}heeler-{F}eynman
  electrodynamics.
\newblock {\em J. Comput. Appl. Math.}, 236(13):3319--3337, 2012.

\bibitem[dlLS19]{MR3900818}
Rafael de~la Llave and Yannick Sire.
\newblock An a posteriori {KAM} theorem for whiskered tori in {H}amiltonian
  partial differential equations with applications to some ill-posed equations.
\newblock {\em Arch. Ration. Mech. Anal.}, 231(2):971--1044, 2019.

\bibitem[Dri60]{MR2613114}
Rodney~David Driver.
\newblock {\em D{ELAY}-{DIFFERENTIAL} {EQUATIONS} {AND} {AN} {APPLICATION} {TO}
  {A} {TWO}-{BODY} {PROBLEM} {OF} {CLASSICAL} {ELECTRODYNAMICS}}.
\newblock ProQuest LLC, Ann Arbor, MI, 1960.
\newblock Thesis (Ph.D.)--University of Minnesota.

\bibitem[Dri63a]{MR0146486}
Rodney~D. Driver.
\newblock A functional-differential system of neutral type arising in a
  two-body problem of classical electrodynamics.
\newblock In {\em Internat. {S}ympos. {N}onlinear {D}ifferential {E}quations
  and {N}onlinear {M}echanics}, pages 474--484. Academic Press, New York, 1963.

\bibitem[Dri63b]{MR0151110}
Rodney~D. Driver.
\newblock A two-body problem of classical electrodynamics: the one-dimensional
  case.
\newblock {\em Ann. Physics}, 21:122--142, 1963.

\bibitem[Dri69]{MR240421}
Rodney~D. Driver.
\newblock A ``backwards'' two-body problem of classical relativistic
  electrodynamics.
\newblock {\em Phys. Rev. (2)}, 178:2051--2057, 1969.

\bibitem[Dri79a]{MR527585}
R.~D. Driver.
\newblock Can the future influence the present?
\newblock {\em Phys. Rev. D (3)}, 19(4):1098--1107, 1979.

\bibitem[Dri79b]{MR553005}
R.~D. Driver.
\newblock Erratum: ``{C}an the future influence the present?'' ({P}hys. {R}ev.
  {D} (3) {\bf 19} (1979), no. 4, 1098--1107).
\newblock {\em Phys. Rev. D (3)}, 20(10):2639, 1979.

\bibitem[dSDL15]{MR3389782}
Daniel~C\^{a}mara de~Souza and Jayme De~Luca.
\newblock Solutions of the {W}heeler-{F}eynman equations with discontinuous
  velocities.
\newblock {\em Chaos}, 25(1):013102, 10, 2015.

\bibitem[DvGVLW95]{MR1345150}
Odo Diekmann, Stephan~A. van Gils, Sjoerd~M. Verduyn~Lunel, and Hans-Otto
  Walther.
\newblock {\em Delay equations}, volume 110 of {\em Applied Mathematical
  Sciences}.
\newblock Springer-Verlag, New York, 1995.
\newblock Functional, complex, and nonlinear analysis.

\bibitem[GAL20]{Mem}
Carlos Garc\'{\i}a-Azpeitia and Jean-Philippe Lessard.
\newblock Free vibrations in a wave equation modeling {MEMS}.
\newblock {\em SIAM J. Appl. Dyn. Syst.}, 19(4):2749--2782, 2020.

\bibitem[GS18]{GOMEZ_PDESurvey}
J.~G{\'o}mez-Serrano.
\newblock Computer-assisted proofs in {PDE}: a survey.
\newblock {\em SeMA Journal}, pages 1--26, 2018.

\bibitem[HCF{\etalchar{+}}16]{MR3467671}
\`Alex Haro, Marta Canadell, Jordi-Llu\'{\i}s Figueras, Alejandro Luque, and
  Josep-Maria Mondelo.
\newblock {\em The parameterization method for invariant manifolds}, volume 195
  of {\em Applied Mathematical Sciences}.
\newblock Springer, [Cham], 2016.
\newblock From rigorous results to effective computations.

\bibitem[HD90]{MR1065250}
Jeffrey~T. Hoag and R.~D. Driver.
\newblock A delayed-advanced model for the electrodynamics two-body problem.
\newblock {\em Nonlinear Anal.}, 15(2):165--184, 1990.

\bibitem[HdlL06a]{MR2240743}
\`A. Haro and R.~de~la Llave.
\newblock A parameterization method for the computation of invariant tori and
  their whiskers in quasi-periodic maps: numerical algorithms.
\newblock {\em Discrete Contin. Dyn. Syst. Ser. B}, 6(6):1261--1300, 2006.

\bibitem[HdlL06b]{MR2289544}
A.~Haro and R.~de~la Llave.
\newblock A parameterization method for the computation of invariant tori and
  their whiskers in quasi-periodic maps: rigorous results.
\newblock {\em J. Differential Equations}, 228(2):530--579, 2006.

\bibitem[HdlL07]{MR2299977}
A.~Haro and R.~de~la Llave.
\newblock A parameterization method for the computation of invariant tori and
  their whiskers in quasi-periodic maps: explorations and mechanisms for the
  breakdown of hyperbolicity.
\newblock {\em SIAM J. Appl. Dyn. Syst.}, 6(1):142--207, 2007.

\bibitem[HdlL16]{MR3501842}
Xiaolong He and Rafael de~la Llave.
\newblock Construction of quasi-periodic solutions of state-dependent delay
  differential equations by the parameterization method {II}: {A}nalytic case.
\newblock {\em J. Differential Equations}, 261(3):2068--2108, 2016.

\bibitem[HdlL17]{MR3736145}
Xiaolong He and Rafael de~la Llave.
\newblock Construction of quasi-periodic solutions of state-dependent delay
  differential equations by the parameterization method {I}: {F}initely
  differentiable, hyperbolic case.
\newblock {\em J. Dynam. Differential Equations}, 29(4):1503--1517, 2017.

\bibitem[Hen21]{MR4292532}
Olivier Henot.
\newblock On polynomial forms of nonlinear functional differential equations.
\newblock {\em J. Comput. Dyn.}, 8(3):309--323, 2021.

\bibitem[HKWW06]{MR2457636}
Ferenc Hartung, Tibor Krisztin, Hans-Otto Walther, and Jianhong Wu.
\newblock Functional differential equations with state-dependent delays: theory
  and applications.
\newblock In {\em Handbook of differential equations: ordinary differential
  equations. {V}ol. {III}}, Handb. Differ. Equ., pages 435--545.
  Elsevier/North-Holland, Amsterdam, 2006.

\bibitem[HLMJ16]{AJPJ}
Allan Hungria, Jean-Philippe Lessard, and J.~D. Mireles~James.
\newblock Rigorous numerics for analytic solutions of differential equations:
  the radii polynomial approach.
\newblock {\em Math. Comp.}, 85(299):1427--1459, 2016.

\bibitem[Joh17]{Johansson2017arb}
F.~Johansson.
\newblock Arb: efficient arbitrary-precision midpoint-radius interval
  arithmetic.
\newblock {\em IEEE Transactions on Computers}, 66:1281--1292, 2017.

\bibitem[KMWZ21]{MR4283203}
Tomasz Kapela, Marian Mrozek, Daniel Wilczak, and Piotr Zgliczy\'{n}ski.
\newblock C{APD}::{D}yn{S}ys: a flexible {C}++ toolbox for rigorous numerical
  analysis of dynamical systems.
\newblock {\em Commun. Nonlinear Sci. Numer. Simul.}, 101:Paper No. 105578, 26,
  2021.

\bibitem[KSW96]{MR1420838}
Hans Koch, Alain Schenkel, and Peter Wittwer.
\newblock Computer-assisted proofs in analysis and programming in logic: a case
  study.
\newblock {\em SIAM Rev.}, 38(4):565--604, 1996.

\bibitem[KW17]{MR3748505}
Tibor Krisztin and Hans-Otto Walther.
\newblock Smoothness issues in differential equations with state-dependent
  delay.
\newblock {\em Rend. Istit. Mat. Univ. Trieste}, 49:95--112, 2017.

\bibitem[Les18]{Lessard2018}
Jean-Philippe Lessard.
\newblock Computing discrete convolutions with verified accuracy via {B}anach
  algebras and the {FFT}.
\newblock {\em Appl. Math.}, 63(3):219--235, 2018.

\bibitem[LMJR16]{MR3545977}
Jean-Philippe Lessard, J.~D. Mireles~James, and Julian Ransford.
\newblock Automatic differentiation for {F}ourier series and the radii
  polynomial approach.
\newblock {\em Phys. D}, 334:174--186, 2016.

\bibitem[LR14]{JPC}
Jean-Philippe Lessard and Christian Reinhardt.
\newblock Rigorous numerics for nonlinear differential equations using
  {C}hebyshev series.
\newblock {\em SIAM J. Numer. Anal.}, 52(1):1--22, 2014.

\bibitem[MPN14]{MR3229655}
John Mallet-Paret and Roger~D. Nussbaum.
\newblock Analyticity and nonanalyticity of solutions of delay-differential
  equations.
\newblock {\em SIAM J. Math. Anal.}, 46(4):2468--2500, 2014.

\bibitem[MPN19]{MR3992063}
John Mallet-Paret and Roger~D. Nussbaum.
\newblock Intricate structure of the analyticity set for solutions of a class
  of integral equations.
\newblock {\em J. Dynam. Differential Equations}, 31(3):1045--1077, 2019.

\bibitem[Nak01]{NAKAO_VerifiedPDE}
M.~T. Nakao.
\newblock Numerical verification methods for solutions of ordinary and partial
  differential equations.
\newblock {\em Numerical Functional Analysis and Optimization},
  22(3-4):321--356, 2001.

\bibitem[NPW19]{MR3971222}
Mitsuhiro~T. Nakao, Michael Plum, and Yoshitaka Watanabe.
\newblock {\em Numerical verification methods and computer-assisted proofs for
  partial differential equations}, volume~53 of {\em Springer Series in
  Computational Mathematics}.
\newblock Springer, Singapore, [2019] \copyright 2019.

\bibitem[Rum99]{Ru99a}
Siegfried~M. Rump.
\newblock {INTLAB - INTerval LABoratory}.
\newblock In Tibor Csendes, editor, {\em {Developments~in~Reliable Computing}},
  pages 77--104. Kluwer Academic Publishers, Dordrecht, 1999.

\bibitem[Rum01]{Rump2001}
Siegfried~M. Rump.
\newblock Computational error bounds for multiple or nearly multiple
  eigenvalues.
\newblock volume 324, pages 209--226. 2001.
\newblock Special issue on linear algebra in self-validating methods.

\bibitem[Rum10]{MR2652784}
Siegfried~M. Rump.
\newblock Verification methods: rigorous results using floating-point
  arithmetic.
\newblock {\em Acta Numer.}, 19:287--449, 2010.

\bibitem[Rum18]{Rump2018}
Siegfried~M. Rump.
\newblock Mathematically rigorous global optimization in floating-point
  arithmetic.
\newblock {\em Optim. Methods Softw.}, 33(4-6):771--798, 2018.

\bibitem[Tre13]{MR3012510}
Lloyd~N. Trefethen.
\newblock {\em Approximation theory and approximation practice}.
\newblock Society for Industrial and Applied Mathematics (SIAM), Philadelphia,
  PA, 2013.

\bibitem[Tuc11]{TUCKER_ValidatedIntroduction}
W.~Tucker.
\newblock {\em Validated numerics: a short introduction to rigorous
  computations}.
\newblock Princeton University Press, 2011.

\bibitem[vdBL15]{VANDENBERG_Dynamics}
J.~B. van~den Berg and J.~P. Lessard.
\newblock Rigorous numerics in dynamics.
\newblock {\em Notices of the AMS}, 62(9):1057--1061, 2015.

\bibitem[vdBL18]{MR3822720}
Jan~Bouwe van~den Berg and Jean-Philippe Lessard, editors.
\newblock {\em Rigorous numerics in dynamics}, volume~74 of {\em Proceedings of
  Symposia in Applied Mathematics}. American Mathematical Society, Providence,
  RI, 2018.
\newblock AMS Short Course: Rigorous Numerics in Dynamics, January 4--5, 2016,
  Seattle, Washington.

\bibitem[vdBS21]{MR4292534}
Jan~Bouwe van~den Berg and Ray Sheombarsing.
\newblock Rigorous numerics for odes using {C}hebyshev series and domain
  decomposition.
\newblock {\em J. Comput. Dyn.}, 8(3):353--401, 2021.

\bibitem[Wal03a]{MR2084753}
Hans-Otto Walther.
\newblock Differentiable semiflows for differential equations with
  state-dependent delays.
\newblock {\em Univ. Iagel. Acta Math.}, (41):57--66, 2003.

\bibitem[Wal03b]{MR2019242}
Hans-Otto Walther.
\newblock The solution manifold and {$C^1$}-smoothness for differential
  equations with state-dependent delay.
\newblock {\em J. Differential Equations}, 195(1):46--65, 2003.

\bibitem[Wal16]{MR3642771}
Hans-Otto Walther.
\newblock Semiflows for differential equations with locally bounded delay on
  solution manifolds in the space {$C^1((-\infty,0],\Bbb R^n)$}.
\newblock {\em Topol. Methods Nonlinear Anal.}, 48(2):507--537, 2016.

\bibitem[Wal21]{MR4261206}
Hans-Otto Walther.
\newblock Solution manifolds which are almost graphs.
\newblock {\em J. Differential Equations}, 293:226--248, 2021.

\bibitem[WdlL20]{MR4120821}
Fenfen Wang and Rafael de~la Llave.
\newblock Response solutions to quasi-periodically forced systems, even to
  possibly ill-posed {PDE}s, with strong dissipation and any frequency vectors.
\newblock {\em SIAM J. Math. Anal.}, 52(4):3149--3191, 2020.

\bibitem[WF49]{MR0032447}
John~Archibald Wheeler and Richard~Phillips Feynman.
\newblock Classical electrodynamics in terms of direct inter-particle action.
\newblock {\em Rev. Modern Physics}, 21:425--433, 1949.

\bibitem[YGdlL]{Per}
Jiaqi Yang, Joan Gimeno, and Rafael de~la Llave.
\newblock Persistence and smooth dependence on parameters of periodic orbits in
  functional differential equations close to an ode or an evolutionary pde.
\newblock https://arxiv.org/abs/2103.05203.

\bibitem[YGdlL21]{MR4287353}
Jiaqi Yang, Joan Gimeno, and Rafael de~la Llave.
\newblock Parameterization method for state-dependent delay perturbation of an
  ordinary differential equation.
\newblock {\em SIAM J. Math. Anal.}, 53(4):4031--4067, 2021.

\end{thebibliography}
\end{document}